\documentclass[final]{siamltex}

\ifx\pdfoutput\undefined
   \usepackage[dvips]{graphicx}
   \else
   \usepackage[pdftex]{graphicx}
   \fi
   \usepackage{epstopdf}

\usepackage{fixltx2e}

\usepackage{graphicx}
\usepackage{amssymb, amsmath}
\usepackage{color}
\usepackage{psfrag}
\usepackage{hyperref}

\def\eps{\varepsilon}

\newcommand\eJGold[1]{#1}
\newcommand\eJG[1]{#1}
\newcommand\ePM[1]{#1}
\newcommand\ePMnew[1]{#1}
\newcommand\ePMnewnew[1]{#1}

\title{Unfoldings of Singular Hopf Bifurcation\thanks{We thank Hinke Osinga, Bernd Krauskopf and Mathieu Desroches for valuable discussions about homoclinic orbit and inclination flips, as well as visualization methods for invariant manifolds. \ePMnew{We moreover thank Bernd Krauskopf and Christian Kuehn for helpful comments on earlier versions of this paper. } This research was partially supported by grants from the National Science Foundation and the Department of Energy.}}

\author{John Guckenheimer, Philipp Meerkamp\thanks{Mathematics Department, Cornell University, Ithaca, NY 1485}}

\begin{document}

\maketitle

\begin{abstract}
Singular Hopf bifurcation occurs in generic families of vector-fields with two slow variables and one fast variable. Normal forms for this bifurcation depend upon several parameters, and the dynamics displayed by the normal forms is intricate. This paper analyzes a normal form for this bifurcation. It presents extensive diagrams of bifurcations of equilibrium points and periodic orbits that are close to singular Hopf bifurcation. In addition, parameters are determined where there is a tangency between invariant manifolds that are important in the appearance of mixed-mode oscillations in systems near singular Hopf bifurcation. One parameter of the normal form is identified as the primary bifurcation parameter, and the paper presents a catalog of bifurcation sequences that occur as the primary bifurcation parameter is varied. These results are applied to estimate the parameters for the onset of mixed-mode oscillations in a model of chemical oscillations.
\end{abstract}





\section{Introduction}

This paper presents an analysis of \emph{singular Hopf bifurcation} in slow-fast dynamical systems with two slow variables. The setting of the problem is as follows. Slow-fast dynamical systems are written
\begin{equation*} 
\begin{split}
\eps \dot{x}
& = 
f(x,y)  \\
\dot{y} & = 
g(x,y), \\
\end{split}
\label{sf}
\end{equation*}
where \ePMnew{the overdot denotes a derivative with respect to a time variable $t$,} $x\in R^m$ are the fast variables, $y\in R^n$ are the slow variables, and $\eps>0$ is a small parameter that is the ratio of time scales. Throughout this paper we assume that there are two slow variables and one fast variable: $m=1$ and $n=2$.

The singular limit of a slow-fast system is a compound entity involving both the slow and fast time scales. Equation~\eqref{sf} becomes a differential algebraic equation upon setting $\eps = 0$. This defines a vector field, the \emph{reduced system}, on the two dimensional \emph{critical manifold} $C$ given by $f=0$ at points where the projection of $C$ onto the space of (slow) variables $y$ is regular. A general procedure that rescales time can be used to extend the reduced system to the \emph{fold} points of $C$. For the desingularized reduced system, there are two types of equilibrium points: ones that satisfy both $f=0$ and $g=0$, and ones that lie on folds but do not satisfy $g=0$. The latter are called \emph{folded singularities}. 

The presence of multiple time scales yields new bifurcation phenomena that continue to be studied extensively. \emph{Singular Hopf bifurcation} in the ``full'' system \eqref{sf} with $\eps > 0$ and \emph{folded saddle-node bifurcation of type II} \cite{SW_canards_in_r3} in the desingularized reduced system together comprise one of these phenomena. Folded saddle-nodes of type II occur at equilibrium points of the desingularized slow flow that satisfy both criteria for a folded singularity: they lie on folds and $g=0$~\cite{SW_canards_in_r3}. In a one parameter family of systems, if an equilibrium point of the desingularized system crosses the folds of $C$, there is a folded saddle-node of type II where it does so. In the full system, there is typically a Hopf bifurcation point at parameter values and state space location that are within $O(\eps)$ of the folded saddle-node of type II. This Hopf bifurcation is singular: its imaginary eigenvalues have a magnitude that is intermediate between the fast and slow time scales. Thus, the folded saddle-node of type II can be viewed as the singular limit of a singular Hopf bifurcation. 

\ePMnew{Singular Hopf bifurcations has been observed in systems modeling action potentials in neurons as well as chemical reactions: examples are a reduced Hogdkin-Huxley model (see \cite{mmo_paper} and \cite{num_cont_canard_seg_slow_fast} for recent results) that has a subcritical singular Hopf bifurcation, a model for chemical reactors introduced by Koper  \cite{koper_paper} \cite{KrupaPopovicKopell} \cite{mmo_paper} that has a supercritical singular Hopf bifurcation, and a model for an autocatalator first studied by Petrov, Scott and Showalter (see \cite{petrov_scott_showalter}, as well as \cite{milik_szmolyan},  \cite{singular_hopf_paper} and \cite{Guckenheimer_Scheper} for more recent studies). All of these models exhibit a type of mixed mode oscillation that is directly related to singular Hopf bifurcation. } \ePMnewnew{Figures in Vo et al~\cite{VoEtAl} moreover suggest that singular Hopf type oscillations are also present in a model of electrical activity and $\mathrm{Ca}^{++}$ dynamics in a pituitary lactotroph, introduced by Toporikova et al ~\cite{ToporikovaEtAl}. }

Guckenheimer~\cite{singular_hopf_paper} presented a normal form for singular Hopf bifurcation and displayed some of the dynamical phenomena that occur in its unfolding. The normal form is 
\begin{equation} 
\begin{split}
\dot{x}
& = 
(y-x^2)/\eps  \\
\dot{y} & = 
z-x \\
\dot{z}& =
-\mu-a x -b y -c z\\
\end{split}
\label{shnf}
\end{equation}
which depends upon the four parameters $\mu, a, b, c$ as well as $\eps$. New dynamical phenomena, not predicted by the Hopf theory in single time-scale systems, are found in this family of vector fields. This is already true in systems with one slow and one fast variable where there is a canard explosion~\cite{benoit}, but the dynamics are far more complicated in systems with two slow variables and one fast variable. Guckenheimer~\cite{singular_hopf_paper} initiated the study of these phenomena. He observed bifurcations of equilibrium points and periodic orbits in the normal form~\eqref{shnf}, and showed that a tangency of invariant manifolds is present. This paper extends that work with a comprehensive analysis of some types of bifurcations that appear in the unfolding of the singular Hopf bifurcation normal form~\eqref{shnf}.  \ePM{The focus of the analysis is on phenomena \emph{local} to the origin, i.e. those that occur at a bounded distance from the origin as $\epsilon\rightarrow 0$. In system ~\eqref{shnf}, local objects or events thus approach the origin as $\epsilon \rightarrow 0$, and are likely to persist under perturbations to system ~\eqref{shnf}.}  \eJGold{The bifurcation diagrams and tables presented here are aids for the analysis of any system with two slow variables in which singular Hopf bifurcation occurs. This is illustrated by an example in Section 5.
We describe our results in the framework of geometric singular perturbation theory ~\cite{jones_lecture_notes}.}

The normal form~\eqref{shnf} is invariant under a scaling transformation that eliminates $\eps$ as a parameter: set $(X,Y,Z,T) = (\eps^{-1/2}x,\eps^{-1} y,\eps^{-1/2}z,\eps^{-1/2}t)$ and $(A,B,C)=(\eps^{1/2}a,\eps b,\eps^{1/2}c)$ 
to obtain
\begin{equation} 
\begin{split}
X'
& = 
Y-X^2 \\
Y' & = 
Z-X \\
Z'& =
-\mu-A X -B Y -C Z\\
\end{split}
\label{resc_shnf}
\end{equation}
Most of the analysis in this paper uses the rescaled version of the normal form. We present a detailed catalog of codimension 1 \eJGold{bifurcations occurring in system \eqref{resc_shnf}, indicate where these bifurcations occur, and how these bifurcations change as additional} parameters are varied. Our strategy is to do this in three stages, first examining a line in the parameter space parallel to the $\mu$-axis, then a two dimensional slice of the parameter space that varies $(\mu,A)$ and finally investigating how the bifurcations within this slice change as $(B,C)$ is varied.

We regard $\mu$ as the \emph{primary} bifurcation parameter and seek to determine the sequences of bifurcations of \eqref{resc_shnf} that occur as $\mu$ varies while $(A,B,C)$ remain fixed. This includes period-doubling bifurcations, folds of periodic orbits, Hopf bifurcations, homoclinic bifurcations, and a tangency of invariant manifolds, in which the two dimensional unstable manifold of an equilibrium in the fold region intersects the repelling slow manifold tangentially. Section 2 presents the bifurcations and a sequence of phase portraits for the family with $(A,B,C)=(-0.05, 0.001, 0.1)$. 

Section 3 studies how the bifurcations identified in Section 2 depend upon the parameter $A$ with $(B,C)=(0.001, 0.1)$ remaining fixed. We present a detailed two dimensional $(\mu,A)$ bifurcation diagram, showing a curve along which tangencies of invariant manifolds occur, bifurcations of the periodic orbit born in the singular Hopf bifurcation, bifurcations of equilibrium points occurring close to the origin, and curves of homoclinic bifurcations. Canard explosions occur near the homoclinic bifurcations. We discuss the codimension two bifurcations that occur in this restricted two parameter family. Due to the numerical problems associated with the large separations of time-scales for very small $\eps$ \cite{mmo_paper} \cite{n_cont_dyn_sys}, we computed the bifurcation curves of periodic orbits as well as the tangency curve in the bifurcation diagrams for values of $\eps$ in the range of approximately  $10^{-2}$ and $10^{-4}$. Section 4 describes how this $(\mu,A)$ bifurcation diagram of system \eqref{resc_shnf} changes as $(B,C)$ is varied and identifies codimension 3 bifurcations. We show some topologically inequivalent $(\mu,A)$  diagrams  \eJGold{in this section} and give a larger portfolio in Appendix B. A table lists all the $\mu$ bifurcation sequences that we found.  

Section 5 discusses a modification of the normal form \eqref{shnf} that makes the critical manifold S-shaped. \eJGold{Koper \cite{koper_paper} introduced a model of chemical oscillations that is equivalent to a subfamily of the modified normal form. We investigate changes that occur in the bifurcation diagrams of \eqref{shnf} due to the} modifications of the system, finding that bifurcations of structures located close to the origin persist and are perturbed only slightly. In the modified system, trajectories that jump from the vicinity of the origin approach another sheet of the critical manifold, flow along the associated attracting slow manifold to a fold and then jump back to the attracting slow manifold which comes close to the origin. This sequence of events constitutes a global return mechanism for trajectories to pass repeatedly near the origin. Mixed mode oscillations~\cite{mmo_paper} occur in this setting with small amplitude oscillations located in the vicinity of the origin and large amplitude oscillations that follow the global returns. We use this example to illustrate how 
the table of bifurcation sequences from section 4 can be used to estimate parameters \eJGold{where a tangency between the unstable manifold of the equilibrium and a repelling slow manifold occurs.}

Section 6 summarizes the results in the paper and discusses aspects of our analysis that remain incomplete. We provide a table listing all bifurcation labels used in this paper together with the full names of the bifurcations in appendix A. Appendix B contains a portfolio of $(\mu,A)$ bifurcation diagrams for system ~\eqref{resc_shnf}. Appendix C contains supplementary data on the positions of codimension 3 bifurcations in system ~\eqref{resc_shnf}. Appendix D comments on the numerical methods we used for the computation of bifurcation diagrams. 
\\


\section{Variations of the primary parameter} 

This section studies bifurcations and phase portraits of~\eqref{resc_shnf} that occur as $\mu$ is varied with $(A,B,C)=(-0.05, 0.001, 0.1)$ fixed. The values of these parameters are selected so that a stable equilibrium point, \eJG{denoted $E_f$, undergoes supercritical Hopf bifurcation. The distance from $E_f$ to the origin is $O(\epsilon^{1/2})$. There are four bifurcations that we identify in the one parameter family obtained from varying $\mu$}:
\begin{itemize}
 \item 
Hopf bifurcation,
\item
Tangency of the unstable manifold of an equilibrium with the repelling slow manifold,
\item
Torus bifurcation of the periodic orbit,
\item
Period doubling of the periodic orbit.
\end{itemize}
We present phase portraits that bracket each of the bifurcations to illustrate the changes that occur in the dynamics of this family as $\mu$ increases.

\begin{figure}[hptb]
\begin{center}
\includegraphics[width=0.8\textwidth]{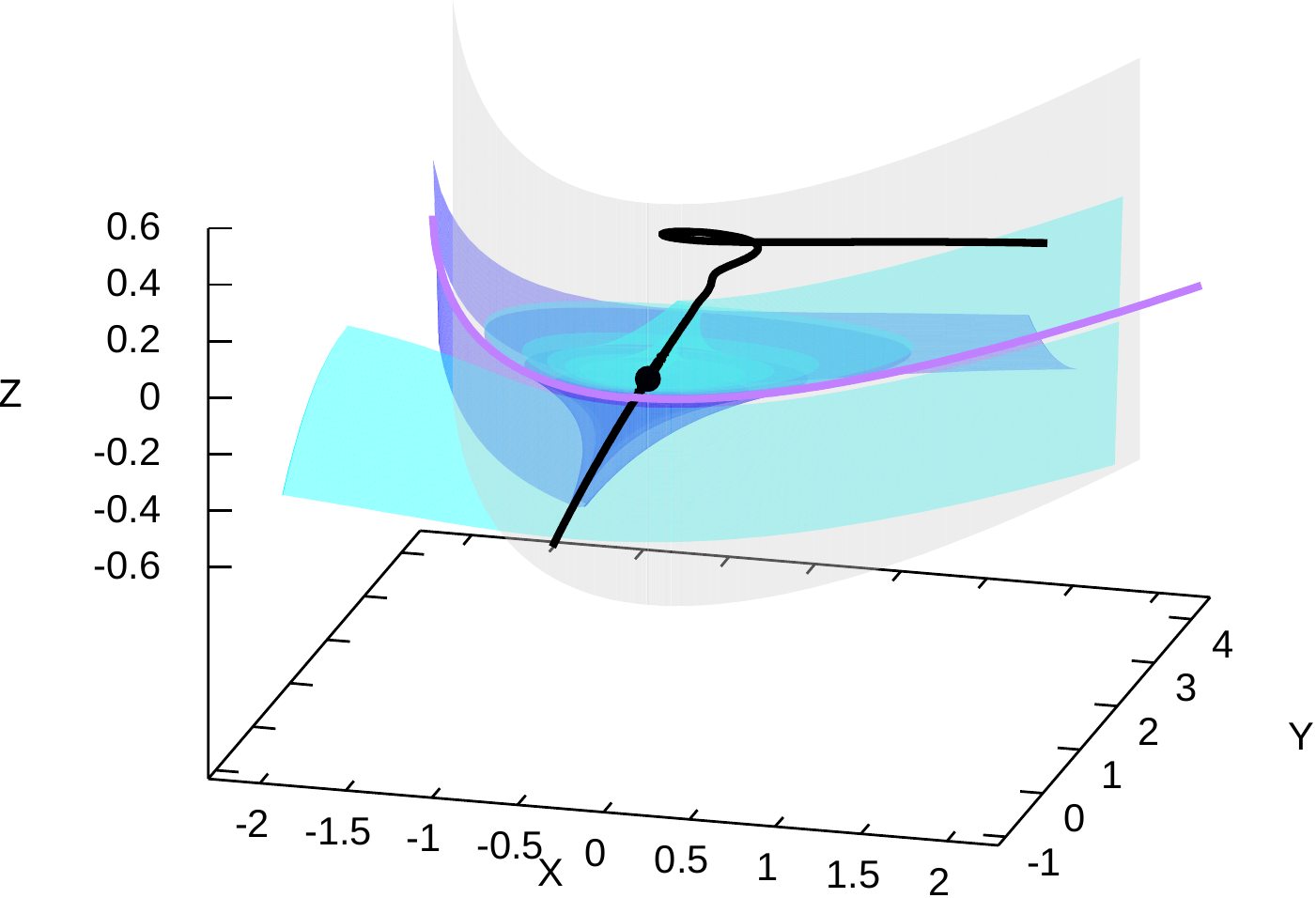} 
\caption{ Phase space when $(\mu,A,B,C)=(0.0, -0.05, 0.001, 0.1)$. The system has an attracting equilibrium with a complex eigenvalue pair at $(X,Y,Z)=(0,0,0)$. Two trajectories tending to the equilibrium along the strong stable eigendirections are drawn with black curves. Parts of $S_r$ that approach the stable manifold of the equilibrium or escape to $X=\infty$ as $t\rightarrow \infty$ are shown in blue. Parts of $S_a$ that approach the equilibrium or tend to $X=-\infty$ as $t\rightarrow \infty$ are shown in cyan. }
\label{fig:after_hopf_mu0Am05B001C1}
\end{center}
\end{figure}

Guckenheimer~\cite{singular_hopf_paper} gives explicit formulas for Hopf bifurcation of~\eqref{resc_shnf} and its first Lyapunov coefficient. Here the Hopf bifurcation occurs for $\mu \approx 0.001246$. Figure~\ref{fig:after_hopf_mu0Am05B001C1} shows the attracting and repelling slow manifolds, denoted $S_a$ and $S_r$ throughout the paper, and the one dimensional strong stable manifold (black) of the equilibrium at a parameter $\mu = 0$ close to the Hopf bifurcation. $S_a$ is drawn in red and magenta, with the top (red) half tending to $E_f$ while the bottom (magenta) half flows to $X = -\infty$. The top half (cyan) of $S_r$ consists of trajectories that flow from $X = -\infty$, while the bottom half (blue) consists of trajectories that originate close to one branch of the \eJG{strong} stable manifold of the equilibrium. The existence of locally invariant slow manifolds of singularly perturbed systems that lie within distance $O(\eps)$ from normally hyperbolic critical manifolds was proved by Fenichel \cite{fenichel_paper}.  \eJGold{Attracting  manifolds have fast foliations consisting of trajectories that approach each other in forward time; repelling manifolds have fast foliations consisting of trajectories that approach each other in backward time.} As $\eps$ tends to zero, the flow on the slow manifolds approaches the flow of the reduced system. Since we have a single fast variable, the slow manifolds are either attracting or repelling. Extensions of these manifolds to the vicinity of the fold curve intersect one another. \ePM{Detailed analyses of such intersections near folded saddles and folded nodes were performed by \eJGold{Benoit~\cite{B}, Szmolyan and Wechselberger ~\cite{SW_canards_in_r3}, Wechselberger ~\cite{Wech}, Guckenheimer and Haiduc~\cite{GHai}, Desroches et al ~\cite{desr_geom_slow_folded_node} and Krupa and Wechselberger ~\cite{KW_local_analysis_folded_saddle_node_singularity}}.} 
%
%
%
%
$S_r$ separates trajectories in the basin of attraction of $E_f$ from trajectories that flow to infinity in forward time.

\begin{figure}[hptb]
\begin{center}
\includegraphics[width=0.8\textwidth]{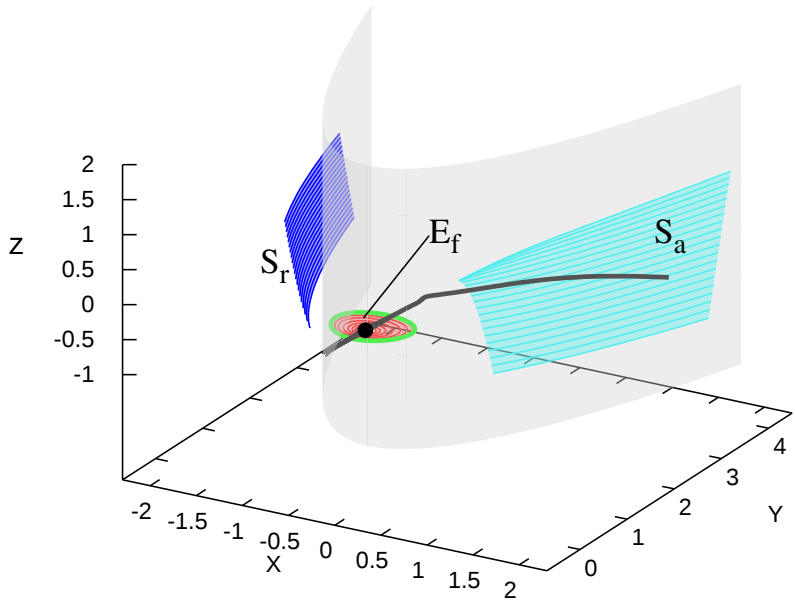}
\psfrag{SR}{$S_r$}
\caption{Phase space when $(\mu,A,B,C)=(0.0012715, -0.05, 0.001, 0.1)$. The critical manifold has the shape of a parabolic sheet and is plotted in gray. Subsets of $S_a$ and $S_r$ are shown in cyan and blue respectively. A black dot marks the position of $E_f$. Its two-dimensional unstable manifold \ePMnew{$W^u(E_f)$} is plotted as a red surface, with a part close to the equilibrium removed. The singular Hopf periodic  orbit is shown as a green curve. The stable manifold  \ePMnew{$W^s(E_f)$} of the equilibrium \ePMnew{$E_f$ is drawn with a black curve, it} leaves the fold region in reverse time: one branch follows $S_r$ and the other remains close to $S_a$ for a while.  }
\label{fig:after_hopf_mu0012715Am05B001C1}
\end{center}
\end{figure}

\begin{figure}[hptb]
\begin{center}
\includegraphics[width=0.8\textwidth]{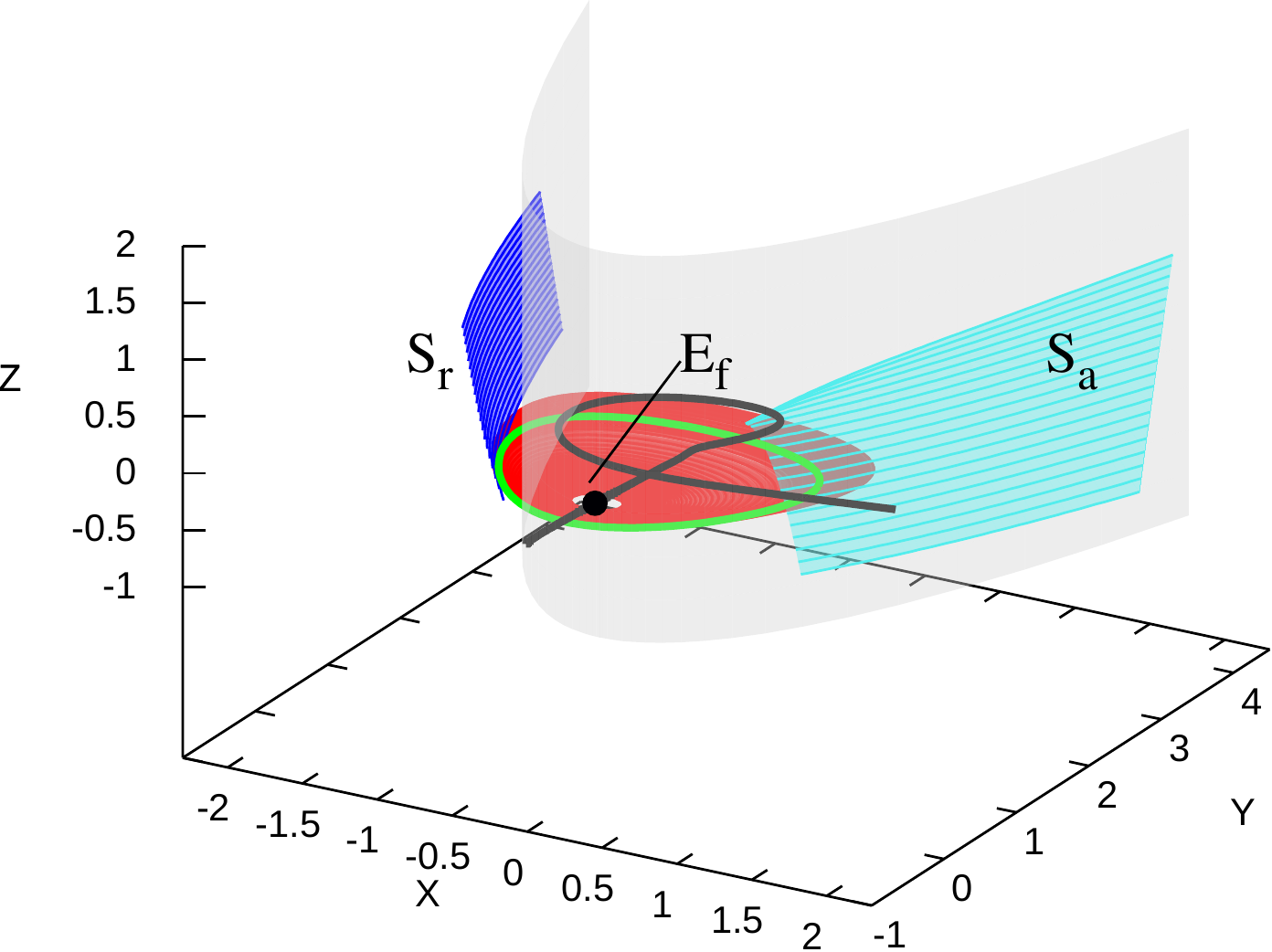} 
\caption{Phase space when $(\mu,A,B,C)=(0.0014975, -0.05, 0.001, 0.1)$, showing the same objects as figure ~\ref{fig:after_hopf_mu0012715Am05B001C1}, with the same color coding. The position of the ``top'' part of \ePMnew{$W^s(E_f)$} changed considerably with respect to figure ~\ref{fig:after_hopf_mu0012715Am05B001C1}. This is to be expected, since the $S_a$ is repelling in reverse time. }
\label{fig:before_tan_mu0014975Am05B001C1}
\end{center}
\end{figure}

\begin{figure}[hptb]
\begin{center}
\includegraphics[width=0.8\textwidth]{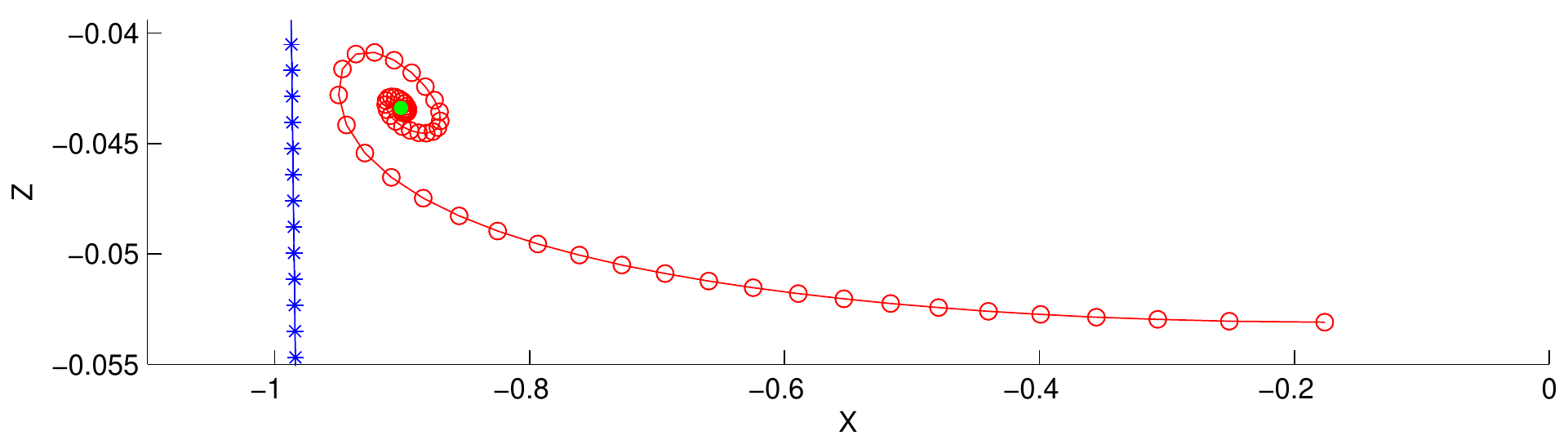} 
\caption{Objects in phase space crossing the plane $Y=0.5$ at $(\mu,A,B,C)=(0.0014975,-0.05,0.001,0.1)$: a selection of trajectories in \ePMnew{$S_r$} are plotted using blue stars, connected with straight lines. The green dot shows the position of the stable periodic orbit's intersection with the plane with $\dot{Y}>0$. The red circles, also connected by straight lines, represent intersections of a single trajectory in \ePMnew{$W^u(E_f)$} with the plane with $\dot{Y}>0$.}
\label{fig:y_section_before_tan}
\end{center}
\end{figure}

Figure~\ref{fig:after_hopf_mu0012715Am05B001C1} displays the phase portrait for $\mu = 0.0012715$. 
The equilibrium point $E_f$ has become a saddle-focus with a two dimensional unstable manifold $W^u(E_f)$ (red) whose boundary is the stable periodic orbit $\Gamma$ (green) born at the Hopf bifurcation. Trajectories in the one-dimensional stable manifold $W^s(E_f)$ are plotted in black. As $\mu$ increases, $\Gamma$ and $W^s(E_f)$ grow in size, getting close to the slow manifolds. At $\mu=0.0014975$, they come close to touching tangentially as illustrated in figures ~\ref{fig:before_tan_mu0014975Am05B001C1} and ~\ref{fig:y_section_before_tan}. The multipliers of $\Gamma$ have become complex, but they still have magnitude smaller than one. Note from figure~\ref{fig:y_section_before_tan} that $W^u(E_f)$ forms a ``scroll'' which gets much closer to $S_r$ than $\Gamma$ does. 

\begin{figure}[hptb]
\begin{center}
\includegraphics[width=0.8\textwidth]{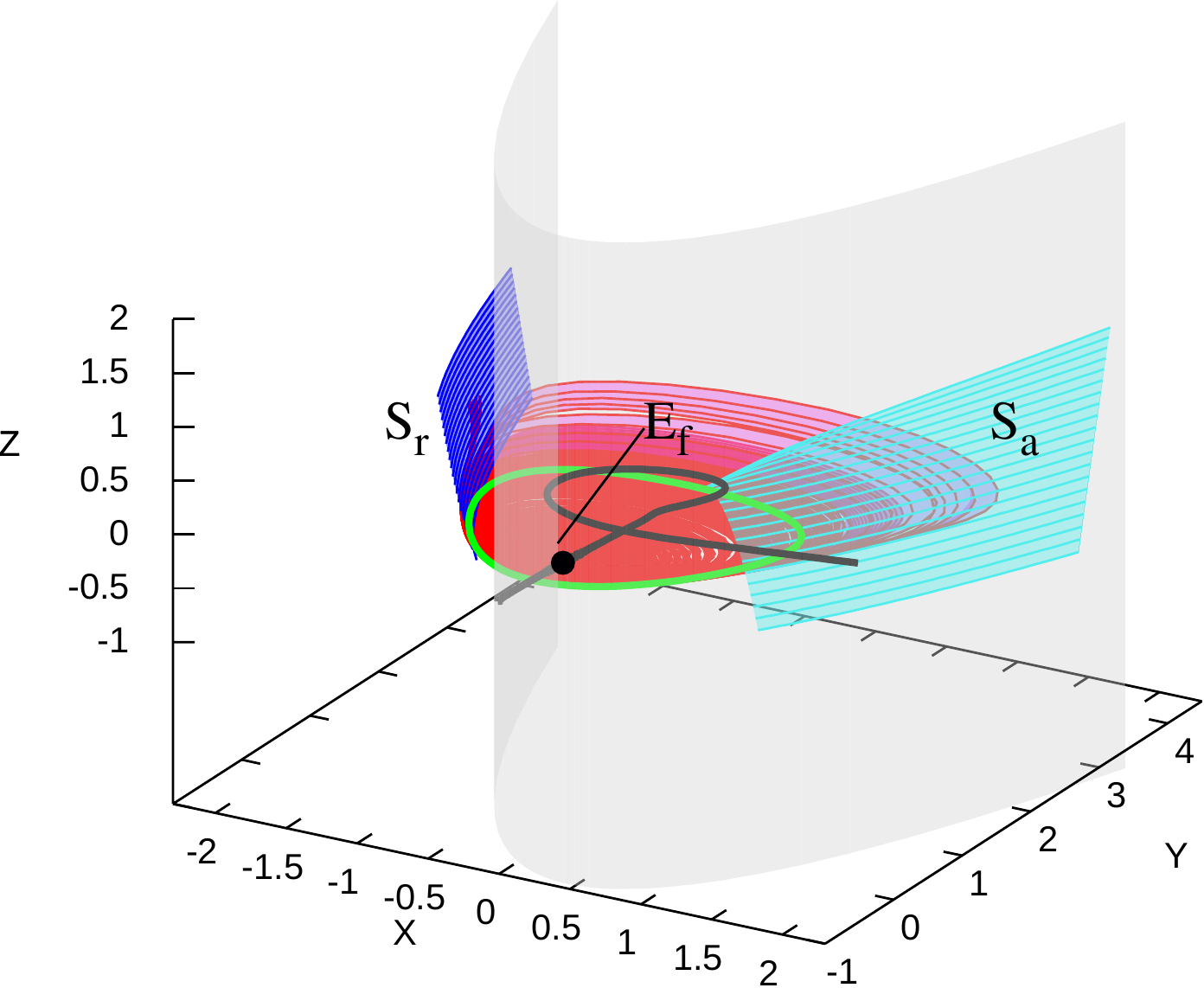}
\caption{Phase space when $(\mu,A,B,C)=(0.0015709, -0.05, 0.001, 0.1)$. In contrast to previous figures, parts of the unstable manifold of the equilibrium \ePMnew{$W^u(E_f)$} near the origin are colored in red and magenta, distinguishing the trajectories that escape to $X=-\infty$ from those that are attracted to the periodic orbit.  }
\label{fig:after_tan_mu001571Am05B001C1}
\end{center}
\end{figure}

$S_r$ and $W^u(E_f)$ begin to intersect transversely at approximately $\mu=0.00156$. This results in bistability for $W^u(E_f)$: some trajectories escape to $X=-\infty$, while others remain in the fold region, and oscillate while approaching $\Gamma$, which still has complex attracting multipliers (see figure ~\ref{fig:after_tan_mu001571Am05B001C1}).

\begin{figure}[hptb]
\begin{center}
\includegraphics[width=0.8\textwidth]{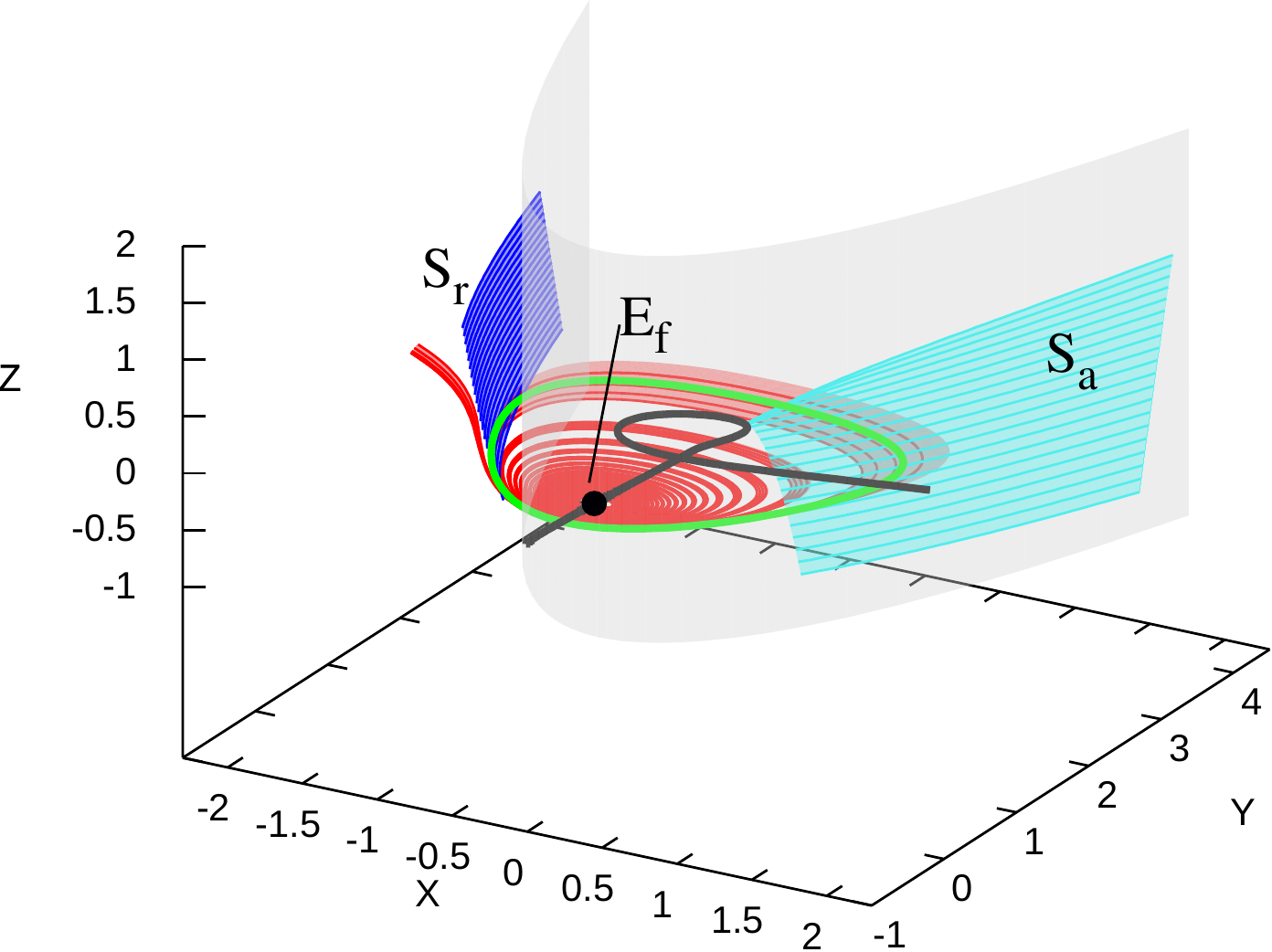}
\caption{Phase space when $(\mu,A,B,C)=(0.0017533, -0.05, 0.001, 0.1)$. The color coding is as in figure ~\ref{fig:after_hopf_mu0012715Am05B001C1}, but here only a part of one fundamental domain of \ePMnew{$W^u(E_f)$} is shown. }
\label{fig:after_tan_mu0017533Am05B001C1}
\end{center}
\end{figure}

Figure ~\ref{fig:after_tan_mu0017533Am05B001C1} shows a phase portrait where $\mu= 0.0017533$. Almost all trajectories in \ePMnew{$W^u(E_f)$} appear to leave the fold region, tending towards $X=-\infty$, with some trajectories following \ePMnew{$S_r$} for extended periods of time on their last turn before leaving the fold region. $W^u(E_f)$ is no longer the boundary of the basin of attraction of $\Gamma$, which remains attracting with complex multipliers.  Note the transverse intersection of $W^u(E_f)$ and $S_r$. 

\begin{figure}[hptb]
\begin{center}
\includegraphics[width=0.8\textwidth]{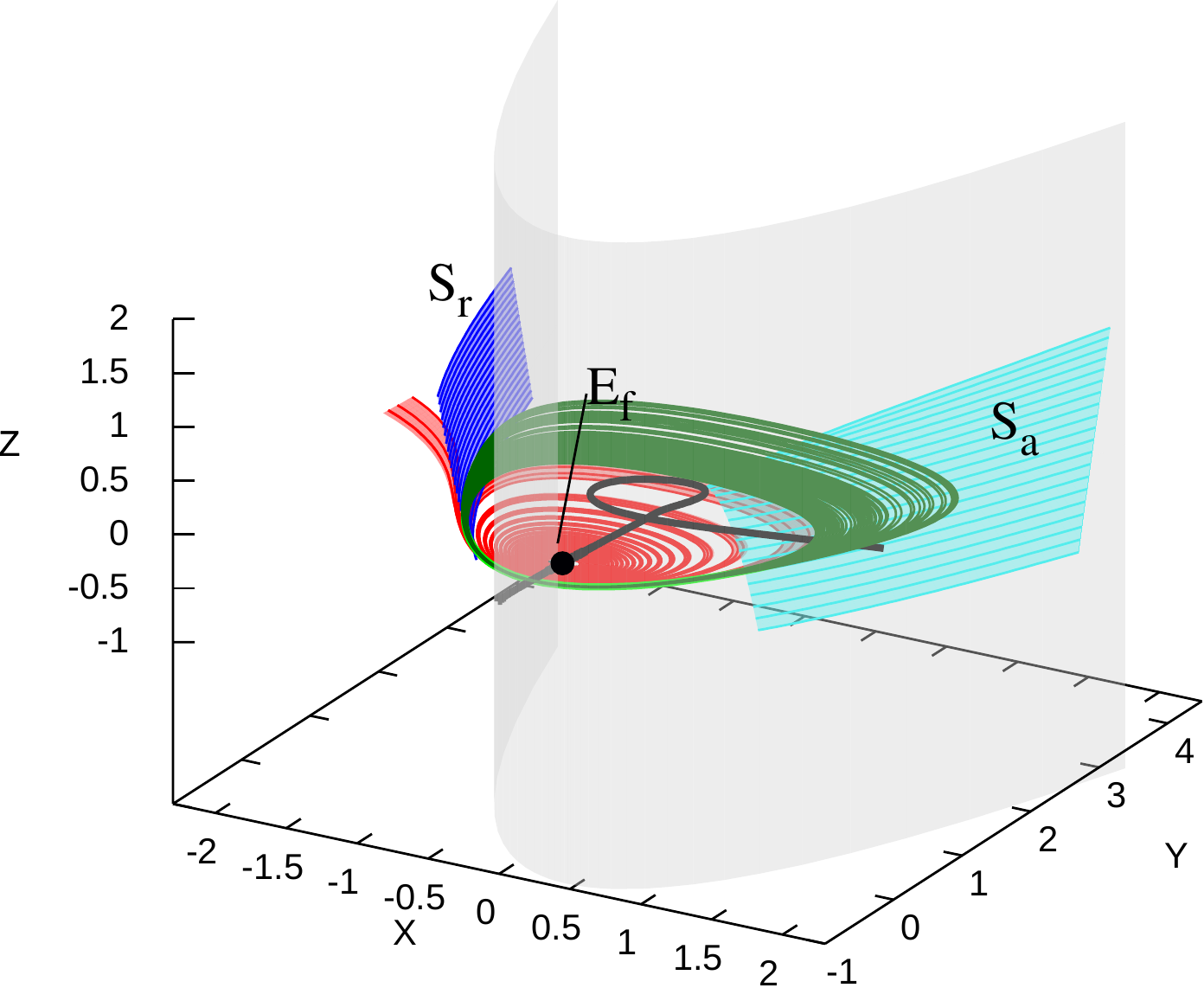}
\caption{Phase space when $(\mu,A,B,C)=(0.0017846, -0.05, 0.001, 0.1)$. A segment of a trajectory on the invariant torus is drawn with a dark-green line. Otherwise, the color coding is as in the previous figure.  \ePM{The periodic orbit is largely obscured by the torus.}}
\label{fig:after_tan_mu0017846087Am05B001C1}
\end{center}
\end{figure}

$\Gamma$ undergoes a torus bifurcation at $\mu=0.0017829$, in which its stability changes from complex attracting to complex repelling, and an invariant torus appears around it (see figure ~\ref{fig:after_tan_mu0017846087Am05B001C1}). 
While numerical computations suggest that at this parameter almost every trajectory in $W^u(E_f)$ tends to $X=-\infty$, the position of the torus and the dynamics close to it may still have an impact on how many ``turns'' trajectories make before exiting the fold region. 
The invariant torus only exists over a short parameter range: at $\mu=0.0017880$, trajectories in the unstable manifold of $\Gamma$ diverge to $X=-\infty$. 

\begin{figure}[hptb]
\begin{center}
\includegraphics[width=0.8\textwidth]{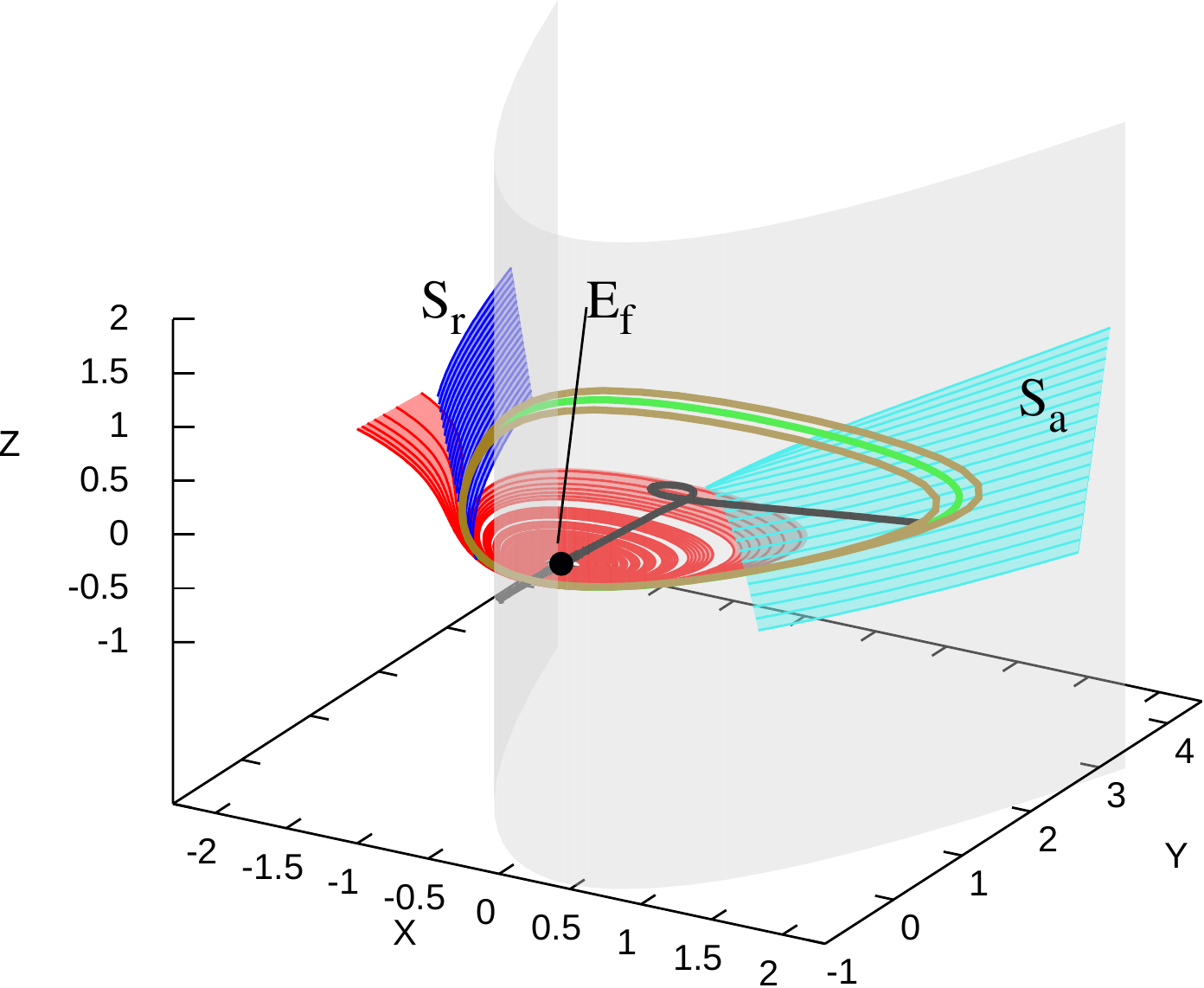}
\caption{Phase space when $(\mu,A,B,C)=(0.0022000, -0.05, 0.001, 0.1)$. A period-doubled periodic orbit is plotted with an olive thick line. The color coding for the other periodic orbit, the slow manifolds, critical manifold, and invariant manifolds of the equilibrium near the origin is as in the earlier figures. A part of the unstable manifold of the equilibrium near the origin not plotted here comes very close to the period-doubled periodic orbit.}
\label{fig:after_tan_mu0022Am05B001C1}
\end{center}
\end{figure}

We note that the multipliers of $\Gamma$ become real repelling for larger values of $\mu$, before $\Gamma$ undergoes a period-doubling bifurcation at approximately $\mu=0.0021910$ (see figure ~\ref{fig:after_tan_mu0022Am05B001C1} for a phase portrait just after the period-doubling).  

\section{Bifurcation structure of a two dimensional parameter slice \texorpdfstring{$(B,C)=(0.001, 0.1)$}{(B,C)=(0.001, 0.1)}  } 

\begin{figure}[htpb]
\begin{center}
\includegraphics[width=1.0\textwidth]{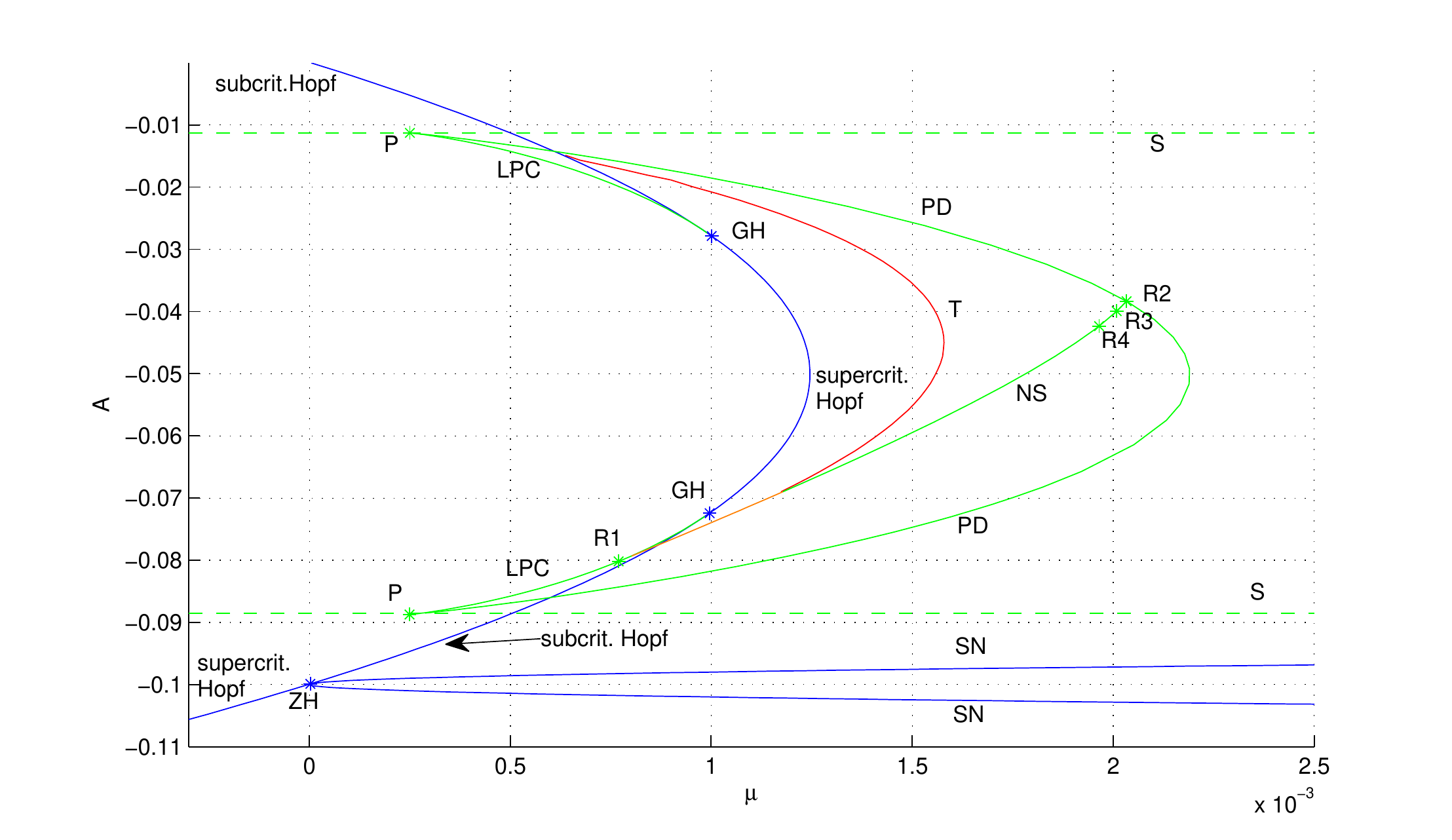} 
\caption{$(\mu,A)$ bifurcation diagram for $(B,C)=(0.001,0.1)$. }
\label{fig:bdiagB001C1_all_aexam}

\end{center}
\end{figure}

This section describes the two-dimensional bifurcation diagram of system ~\eqref{resc_shnf}, with parameters $(B,C)=(0.001, 0.1)$ fixed and $(\mu,A)$ varying. The diagram in figure~\ref{fig:bdiagB001C1_all_aexam} is typical in the sense that all other generically occurring $(\mu,A)$ diagrams share many of its codimension 1 and 2 phenomena and bifurcations. The codimension one bifurcations with the symbols that label them in figure~\ref{fig:bdiagB001C1_all_aexam} are
\begin{itemize}
\item saddle-node bifurcation (SN)
\item (singular) Hopf bifurcation (Hopf)
\item period-doubling bifurcations (PD)
\item fold of periodic orbit bifurcations (LPC)
\item torus bifurcations (NS)
\item tangency of invariant manifolds (T)
\item canard explosions ending in homoclinic bifurcations (S)
\end{itemize}

The codimension two bifurcations are
\begin{itemize}
\item resonances of periodic orbits (R1,R2,R3,R4)
\item degenerate homoclinic bifurcations (P)
\item zero Hopf bifurcations (ZH)
\item generalized Hopf bifurcations (GH)
\end{itemize}
We briefly discuss each of these bifurcations in turn.\\

{\bf Saddle-node bifurcation.} Recall that we assume that $A$ and $C$ are $O(\eps^{1/2})$ and that $B$ is $O(\eps)$. System~\eqref{resc_shnf}  can have up to two equilibria. There is only one equilibrium point near the origin unless $A+C=O(\eps)$, in which case there are two. There is a saddle-node bifurcation (labeled SN in the diagrams) at $\mu =(A+C)^2/(4B)$ and $(X,Y,Z)=(X_e,X_e^2,X_e)$ with $X_e=(A+C)/2B$. The solution curve of $\mu =(A+C)^2/(4B)$ bounds a thin region of parameters where no equilibria exist. \\

{\bf Hopf bifurcation and its codimension 2 degeneracies.} The formulas for Hopf bifurcation loci  in the fold region and its codimension 2 degeneracies are complicated expressions, and are best understood by disregarding terms that are of higher order in $\eps$. 
For all $B\neq 0$, the system has only one Hopf bifurcation that occurs at O(1) distance from the origin, occurring at  $\mu=-A^2/2-AC/2+O(\epsilon^{3/2})$ and $(X,Y,Z)\approx(A/2,A^2/4,A/2)+O(\epsilon^{1/2})$, making it a singular Hopf bifurcation. The Hopf bifurcation is degenerate at a zero Hopf bifurcation (marked ZH in the diagrams), where the Hopf equilibrium point has a zero eigenvalue in addition to the purely imaginary eigenvalue pair, and at two generalized Hopf bifurcations (GH), where the first Lyapunov coefficient vanishes. 
The Hopf center manifold is attracting for values of $A$ greater than that of the zero Hopf bifurcation, and repelling for values of $A$ smaller than that of the zero Hopf bifurcation. 
The criticality of the Hopf bifurcation changes at each of these these codimension 2 degeneracies. Supercritical [subcritical] Hopf bifurcations at values of $A$ greater than that of the ZH bifurcation give rise to periodic orbits when the Hopf curve is crossed with $\mu$ increasing [decreasing], while at smaller values of $A$ supercritical [subcritical] Hopf bifurcations give rise to periodic orbits as the Hopf curve is crossed with $\mu$ decreasing [increasing]. \\

{\bf Period-doubling bifurcations.}  Numerical continuation with MatCont ~\cite{matcont} and AUTO ~\cite{AUTO} show that $\Gamma$ undergoes a period-doubling bifurcation (labeled PD in bifurcation diagrams) as $\mu$ is increased. In the rescaled system~\eqref{resc_shnf}, these bifurcations typically occur on periodic orbits that are $O(1)$ distant from the origin. Near the endpoints of the period-doubling bifurcation curve, \ePMnew{where the PD curve meeds a homoclinic bifurcation,} the amplitude becomes even larger and the bifurcating orbits of system~\eqref{shnf} no longer shrink to the origin as $\eps \to 0$. \ePMnew{This degeneracy is described further in the paragraph on ``Degenerate homoclinic bifurcations'' below. }\\

{\bf Fold of periodic orbits bifurcations.} A fold of periodic orbits bifurcation curve (labeled LPC) emanates from each of the generalized (singular) Hopf bifurcations (GH). This bifurcation, in which two periodic orbits annihilate each other as $\mu$ decreases, is initially local to the origin in system~\eqref{shnf}, \ePM{i.e. remains at bounded distance from the origin as $\epsilon\rightarrow 0$,} but later becomes non-local in a canard explosion, and ends where the \eJGold{LPC curve} meets a homoclinic bifurcation  \ePMnew{(see the paragraph on ``Degenerate homoclinic bifurcations''  for more details)}. \\

{\bf Torus bifurcations.} A locus of torus bifurcations (NS) emanates from the zero Hopf bifurcation (ZH). Associated with this are attracting invariant tori that surround the bifurcating singular Hopf periodic orbit. \eJGold{These tori apparently exist only} in a thin strip of parameters to the right of the torus bifurcation curve. \\

{\bf Tangency of invariant manifolds.} In figure ~\ref{fig:bdiagB001C1_all_aexam}, $E_f$ has a two dimensional unstable manifold in the parameter region bounded below by the saddle-node curve, and on the left by the Hopf curve. In part of this region, $W^u(E_f)$ intersects $S_r$. This intersection is tangential at the tangency curve (T), where the manifolds begin to intersect as $\mu$ increases. Note that there are values of $A$ where the manifolds intersect immediately after the Hopf bifurcation. See appendix D for details on the computation of the tangency curve. \\

\begin{figure}[hptb]
\begin{center}

\includegraphics[width=0.85\textwidth]{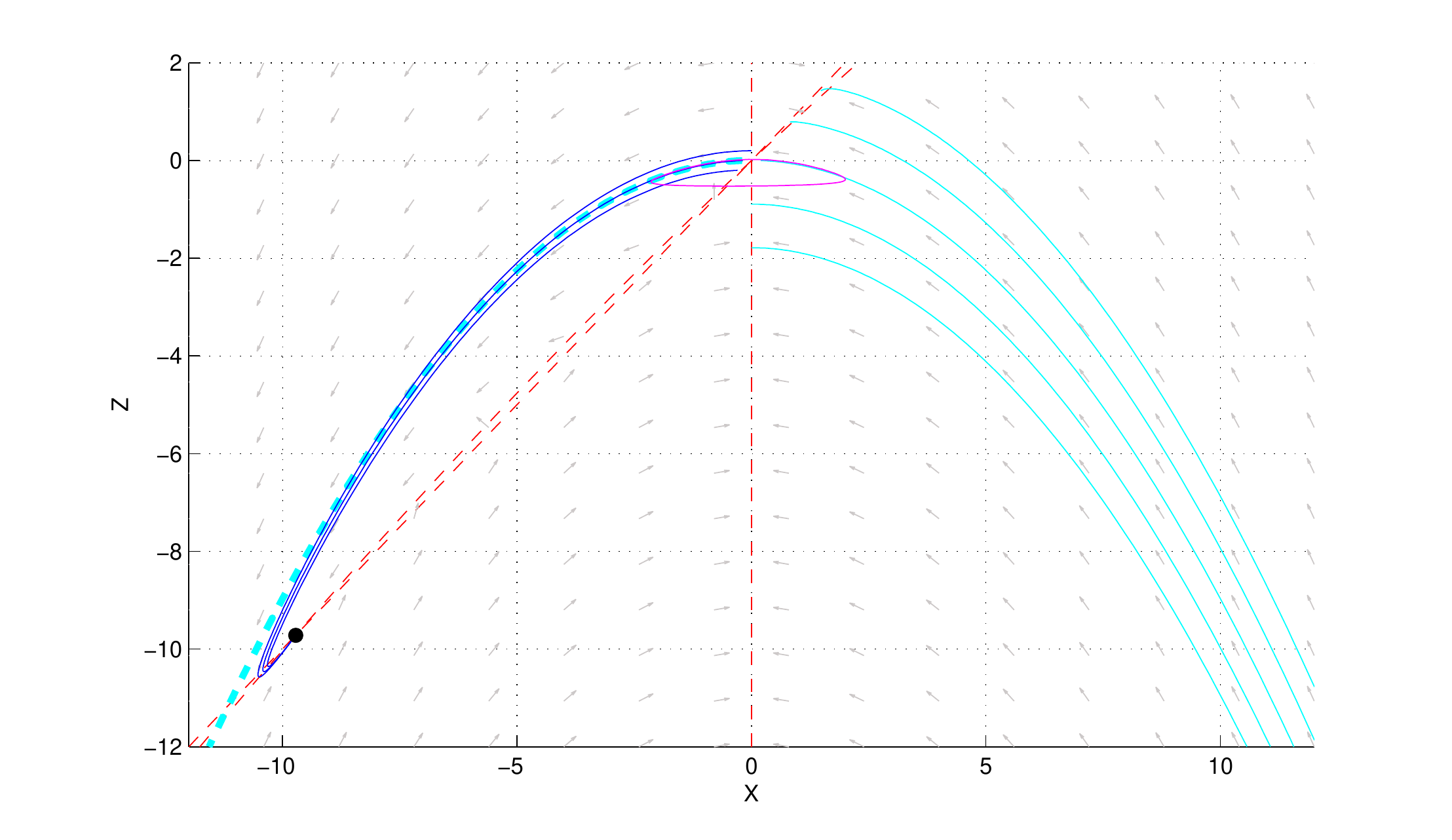} 

\includegraphics[width=0.85\textwidth]{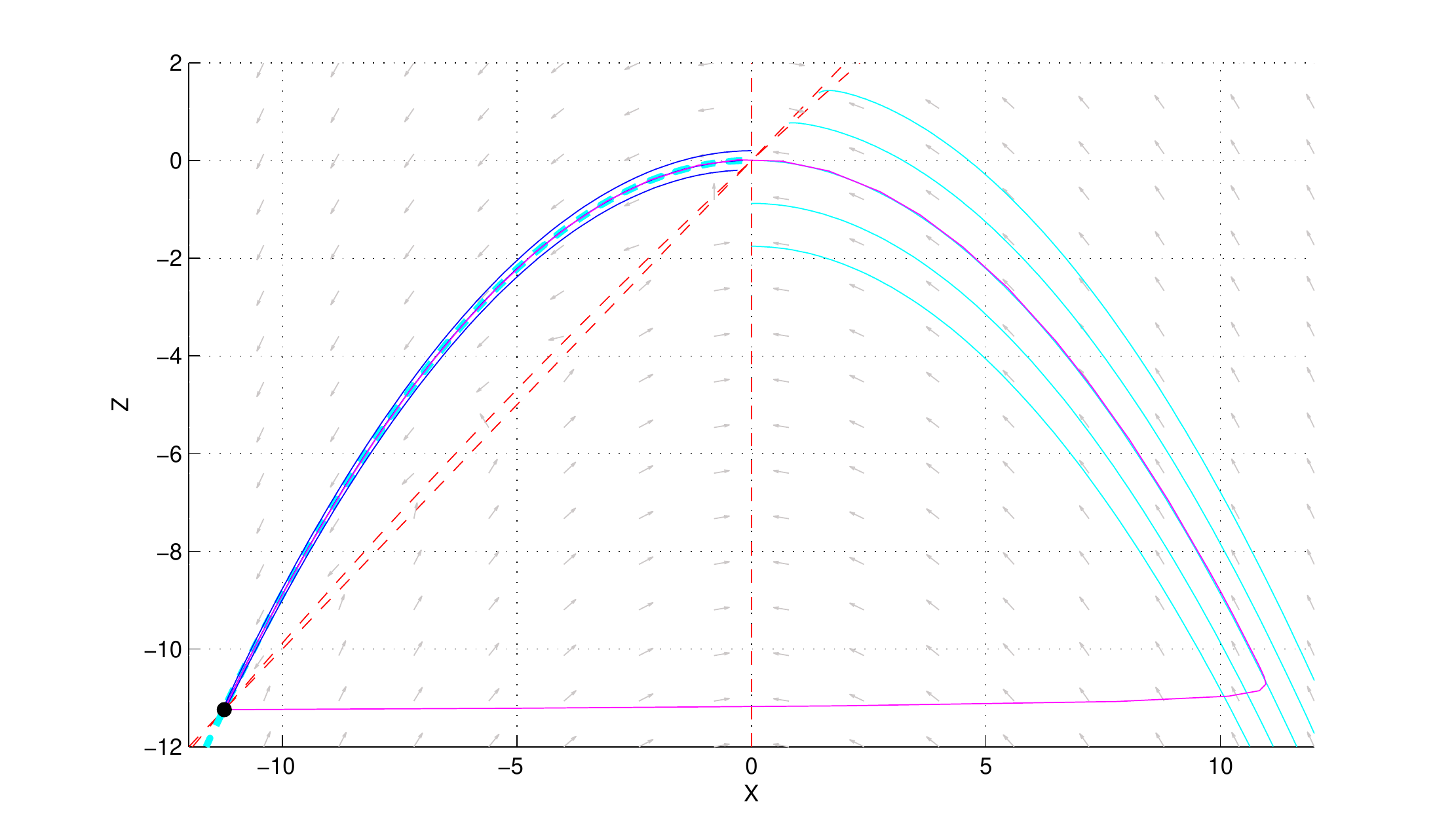} 
\caption{Phase portraits for parameters $(\mu,A,B,C)=(2.5008\times 10^{-6},-0.090281, 0.001,0.1)$ (top) and $(\mu,A,B,C)=(2.5008\times 10^{-6},-0.088758, 0.001,0.1)$ (bottom), 
showing the projection of a periodic orbit (magenta, computed in the full system with MatCont) onto the slow flow.
The nullclines of the slow flow as well as  the fold curve are drawn using red dashed lines. Segments of trajectories of the slow flow are plotted in blue and cyan, the two middle ones ending at the origin. A reflection of the middle cyan trajectory about the $Z$ axis is plotted as a dashed cyan line. The position of \eJGold{the equilibrium far from the origin} on the repelling sheet of the \ePMnew{critical} manifold is marked with a dot. 
}
\label{fig:orb_grow}
\end{center}
\end{figure}

{\bf Canard explosions and homoclinic bifurcations.} The periodic orbit $\Gamma$ grows rapidly in canard explosions. The canard explosions consist of periodic orbits that have long segments lying near $S_r$. The small parameter ranges in which they occur are determined by the geometry of the slow flow as illustrated in figure ~\ref{fig:orb_grow}. The slow flow is computed by differentiating the equation of the critical manifold in system~\eqref{shnf} with respect to \eJG{$t$}, and using the result and the equations for the derivatives of the slow variables to obtain equations for $\dot{x}$ and $\dot{z}$ in $x$ and $z$ only. After rescaling the vector field so that the direction reverses on the repelling sheet of the critical manifold, we obtain the following desingularized slow flow equations:
\begin{equation} 
\begin{split}
x'
& = 
z-x \\
z'& =
\eJGold{-2x\,(\mu+a x + b x^2 +c z)}\\
\end{split}
\label{shnf_slow_flow}
\end{equation}
Trajectories that follow $S_r$ can jump at any location along the manifold. The singular limit of trajectories of system~\eqref{shnf} approach concatenations of trajectories of the slow flow together with segments parallel to the $x$-axis that begin on the repelling sheet of the critical manifold and end on the attracting sheet of the critical manifold. Since the critical manifold is a parabolic cylinder symmetric with respect to reflection in the $(y,z)$ plane, the singular limits of some periodic orbits consist of trajectories of the slow flow that connect points $(x,z)$ to points $(-x,z)$ concatenated with the horizontal segment from $(-x,z)$ to $(x,z)$. We call these periodic orbits with a single fast segment \emph{simple}. \eJG{We refer to their singular limits as singular cycles.}

We expect the simple periodic orbits to grow rapidly where the slow flow equations satisfy $$\frac{\dot{z}(x,z)}{\dot{x}(x,z)}\approx -\frac{\dot{z}(-x,z)}{\dot{x}(-x,z)}.$$
With the simplifying assumption $\mu=0$, substitution of the slow flow equations into the ansatz gives that when $\mu=0$, canard explosions are likely to occur at parameters satisfying 
\begin{equation}
a^2+ac+b=0,
\label{infloc_eqn}
\end{equation}
and that the trajectory segments on the slow manifolds of system~\eqref{resc_shnf} can approximately be parametrized by $(X,Y,Z)=(X,X^2,AX^2)$. Note that $\eps(a^2+ac+b) = A^2+AC+B$ in the rescaling from system~\eqref{shnf} to system~\eqref{resc_shnf}. Numerical calculations with MatCont and AUTO show that these estimates are very accurate, and that the location of periodic orbits is insensitive to the small value of $\mu$. In figure ~\ref{fig:bdiagB001C1_all_aexam}, canard explosions occur along lines that are nearly horizontal, and end in homoclinic bifurcations with equilibria that are far from the origin. 
Dashed green lines labeled S are drawn in figure ~\ref{fig:bdiagB001C1_all_aexam} at $A=\sqrt{C^2-4 B}$ \eJG{as approximations to the region where the slow flow has singular cycles.}   \\

{\bf Resonances of periodic orbits. }There are 1:1, 1:4, 1:3 and 1:2 resonances on the torus bifurcation curve (labeled R1, R4, R3 and R2 respectively), where $\Gamma$ has multipliers on the unit circle with arguments $0,\pi/4,\pi/3$ and $\pi/2$ respectively. At the 1:1 resonance, the torus bifurcation curve ends at a fold of periodic orbits bifurcation locus in a bifurcation that is analogous to the Bogdanov-Takens bifurcation for 2-dimensional vector fields~\cite{Takens,GH}. At the 1:2 resonance, the torus bifurcation curve ends at the period-doubling curve \eJGold{in another bifurcation analogous to the Bogdanov-Takens bifurcation with symmetry~\cite{Takens,GH}.} \\

{\bf Degenerate homoclinic bifurcations. } The homoclinic bifurcation is degenerate at two points, labeled P in the $(\mu,A)$ bifurcation diagrams. At each of these points, a fold of periodic orbits curve (LPC) and a period-doubling curve (PD) meet the curve of homoclinic bifurcations. We believe that these codimension 2 phenomena are due to homoclinic orbit flips (top point) and homoclinic inclination flips (lower point). \ePM{Unfoldings of these bifurcations were developed by Sandstede ~\cite{SAN_thesis} (orbit flip) and Homburg et al ~\cite{HKN_hom_doubl_cascades} (inclination flip) and Kisaka et al.~\cite{KKO_N_homocl},~\cite{KKO_supplement} (inclination flip). We did not investigate the points P in detail.} \\

\section{Catalog of bifurcation diagrams in one and two dimensional parameter slices}

This section first describes curves of codimension 1 bifurcations in $(\mu,A)$ bifurcation diagrams for different fixed values of $(B,C)$. This includes information about the codimension 2 bifurcations at which codimension 1 bifurcation curves meet or end. We then complete our study of the bifurcations of system~\eqref{resc_shnf} by investigating how the $(\mu,A)$ bifurcation diagrams change as $B$ and $C$ vary. These changes occur when we encounter points of codimension larger than two and when we encounter intersections \eJGold{of the loci} of codimension one and two bifurcations that are not transverse. Most of the changes we find are of the latter type. Indeed, we identify here only a single degeneracy of the first type, which occurs when $B=0$.\\

{\bf Equilibrium point bifurcations.} 
Unless $B=0$, system ~\eqref{resc_shnf} has a saddle-node bifurcation at $\mu=(A+C)^2/(4B)$, whose properties are described in section 3. Similarly,  there is only one curve of local Hopf bifurcations for $B\neq 0$, occurring at  $\mu\approx-A^2/2-AC/2+O(\epsilon^{3/2})$. It is degenerate at a zero Hopf (ZH) bifurcation at $A=C(B-1)$. 
A computation of the normal form coefficients using notation from Kuznetsov~\cite{kuznetsov_book} shows that if $B<0$, then $s=1, \, \theta<0$, and the system has homoclinic orbits, invariant tori, and invariant spheres, which are local in the vicinity of the zero Hopf bifurcation. 
If $B>0$, then $s=1, \, \theta>0$, and the dynamics in the vicinity of the zero Hopf bifurcation is simpler. \\

The Hopf bifurcation is also degenerate at generalized Hopf bifurcations (GH), where the first Lyapunov coefficient vanishes. Simplifying the expression for the first Lyapunov coefficient by eliminating terms that are of higher order in $\epsilon$, we observe that there are generalized Hopf bifurcations approximately when  $A^2  + AC+2B=0$ has a real solution, i.e. whenever $C^2-8B\geq 0$. 
We note that equilibrium points at $O(1/\eps)$ distance from the origin can undergo Hopf bifurcation as well. \\

{\bf Periodic orbit bifurcations.} Period-doubling bifurcations (PD) exist in those $(B,C)$ regions where the slow flow has singular cycles. The bifurcation loci begin and end at points labeled P, where the periodic orbit that is to double undergoes a canard explosion. At these points, the period-doubling bifurcation curve meets a fold of periodic orbits curve, in what seems to be a tangential way. The period-doubling occurs as $\mu$ increases through the bifurcation parameter value, and the first period-doubling bifurcation may be followed by further period-doublings. 
Fold of periodic orbits (LPC) curves exist in those $(B,C)$ regions where the slow flow has singular cycles. There is one fold of periodic orbits curve if no generalized Hopf bifurcation is present, and the bifurcation curves begin and end at meeting points labeled P, where they connect with period-doubling bifurcations. If two generalized Hopf bifurcations are present, then the fold of periodic orbits curves begin at generalized Hopf bifurcations, and end at points marked P. 
Given $(B,C)$, the sign of $B$ determines whether the LPC curve lies to the left ($B>0$) or to the right ($B<0$) of the Hopf bifurcation curve. Canard explosions of the slow flow occur for parameters near those that satisfy $A^2+AC+B=0$ where the slow flow has singular cycles.
These canard explosions typically end in homoclinic bifurcations with equilibrium points that are at $O(\epsilon)$ distance from the origin. Being non-local in many cases, we did not investigate homoclinic bifurcations thoroughly and plot homoclinic bifurcation curves only in 
figure~\ref{fig:bdiag_Bm01C1_all_aexam_lpo_w_hom} which gives the $(\mu,A)$ bifurcation diagram for $(B,C)=(-0.01,0.1)$.

Invariant tori exist close to the torus bifurcation curves marked NS. 
If $B<0$, normal form theory predicts that these tori are unstable close to the zero Hopf bifurcation ~\cite{kuznetsov_book}, and exist for values of $A$ greater than that of the torus bifurcation curve. A more detailed analysis of the stability and persistence of the torus bifurcations was not performed. 
When $B<0$, the torus bifurcation curve starts at the zero Hopf bifurcation, and ends at a 1:2 resonance (labeled R2), where it meets the period-doubling curve. The torus bifurcation curve intersects no other bifurcation curves included in the diagrams. The defining equation that characterizes torus bifurcations is also satisfied by ``neutral'' periodic orbits with two multipliers whose product is one. Consider first parameters $(B,C), B>0$ for which there is no generalized Hopf bifurcation. A curve $L_{(B,C)}$ of neutral periodic orbits emanates from the Hopf bifurcation. This curve is shown in figure ~\ref{fig:bdiagB001C065_all_aexam}, but is omitted in all other bifurcation diagrams. Close to the zero Hopf bifurcation, the multipliers are real, and the periodic orbits are saddles. The curve $L_{(B,C)}$ can intersect the period-doubling bifurcation curve in 1:2 resonances or meet a fold of periodic orbits bifurcation curve tangentially in  1:1 resonances where the multipliers become complex. Both events can occur more than once in the same $(\mu,A)$ diagram, due to bending of the curve  $L_{(B,C)}$, see e.g. figure ~\ref{fig:bdiagB001C074_all_aexam} for two 1:2 resonances, or figure ~\ref{fig:bdiagB001C065_all_aexam}. If there is a generalized Hopf bifurcation and $B>0$, then there is typically one curve of torus bifurcations, beginning at a fold of periodic orbits bifurcation at a 1:1 resonance, and ending at  an intersection with a period-doubling bifurcation curve at a 1:2 resonance. 
Independent of whether a generalized Hopf bifurcation is present or not, when $B>0$ and $C>0$, as the absolute value of $C$ decreases towards the parameter $C=C(B)$ where an isola vanishes, the curve $L_{(B,C)}$ tends to move along the families of periodic orbits, from the singular Hopf bifurcation towards  and beyond the period-doubling bifurcation. 
See figures ~\ref{fig:bdiagB001C16_all_aexam} to ~\ref{fig:bdiagB001C065_all_aexam} for a sequence of diagrams illustrating this trend. If $B>0$ and $C<0$, then the families of periodic orbits are traversed away from the Hopf bifurcation as $C$ increases towards the point where the isola vanishes.  \\

{\bf Tangency of invariant manifolds.} We now discuss intersections of $W^u(E_f)$ with $S_r$. 
Regions in parameter space where these manifolds intersect are separated from those where they do not intersect by tangencies of invariant manifolds as well as by three codimension 1 bifurcations: saddle-node bifurcations, Hopf bifurcations, and fold of periodic orbits bifurcations.
The tangency curve is in many cases very close to torus bifurcations or period-doubling bifurcations, to the point of being indistinguishable from these bifurcation curves in some of the $(\mu,A)$ bifurcation diagrams. Nonetheless, the tangency of invariant manifolds is a bifurcation in its own right, and does not coincide with period-doubling or torus bifurcations in system \eqref{resc_shnf} except at isolated points of the $(\mu,A)$ diagrams. 

The tangency curve is situated in the $(\mu,A)$ diagrams in different ways for $B>0$ and $B<0$. 
If $B<0$, $E_f$ leaves an $O(1)$ neighborhood of the origin in system \eqref{resc_shnf} when $A+C=O(\eps)$ with $A+C<0$. The tangency curve begins at the saddle-node curve at a value of $A$ slightly more negative than that of the zero Hopf bifurcation, indistinguishably close in most of the diagrams of this paper. 
In this region of parameter space, $W^u(E_f)$ may approach the one-dimensional stable manifold of \eJGold{the second} equilibrium in the vicinity of the origin and then leave the fold region, or can ``get caught'' in invariant structures around the two equilibrium points involved in the zero Hopf bifurcation, cf. ~\cite{kuznetsov_book} figure 8.22.
%
The tangency curve may intersect the the PD curve, but systematic numerical computations suggest that the tangency curve does not lie visibly to the right of it at the resolution of the $(\mu,A)$ bifurcation diagrams included in this paper. 
The tangency curve ends at the top point $P$, where the boundary of the region with an intersection of invariant manifolds begins to coincide with the fold of periodic orbits bifurcation curve, and the intersection begins to be transverse at the boundary parameters.

If $B>0$, there is a saddle-focus in the vicinity of the origin in the region $R_{(B,C)}$ that lies to the right of the Hopf curve and above the saddle-node curve. If now $C<0$, then \ePMnew{$W^u(E_f)$} and \ePMnew{$S_r$} intersect transversely in the entire region $R_{(B,C)}$. If  $C>0$, then in some parts of  the region $R_{(B,C)}$ the two invariant manifolds do not intersect. More precisely, there is a tangency curve, beginning on the Hopf bifurcation curve, very close to the torus bifurcation curve or at a larger value of $A$, lying to the right of the Hopf curve until it meets the Hopf curve again at a larger value of $A$, very close to the period-doubling curve or at a smaller value of $A$. The tangency curve can intersect the torus and period-doubling curves for large enough values of $C=C(B)$, and there may be interactions between  \ePMnew{$W^u(E_f)$} and the attracting tori for $(\mu,A)$ slightly below the torus bifurcation curve. Systematic numerical computations suggest that the tangency curve does not lie visibly to the right of the period-doubling curve at the resolution of the $(\mu,A)$ bifurcation diagrams included in this paper. It never lies much to the right of the torus bifurcation curve either, see figure ~\ref{fig:bdiagB001C077_all_aexam} for a case in which the tangency curve lies visibly to the right of the torus bifurcation curve. \\

{\bf Remark: local and non-local bifurcations.} The bifurcations described above are, with two exceptions, local to an $O(\eps^{1/2})$ neighborhood of the origin when scaled back to the system~\eqref{shnf}. The first exception is the saddle-node bifurcation (labeled SN), which is local in the vicinity of the zero Hopf bifurcation (ZH), but is  generally non-local. Note also that the range of values of $A$ for which saddle-node bifurcations occur for $(B,C)$ fixed and $\mu = O(\eps)$ is $O(\eps)$, while most other codimension one bifurcations are of interest over an $O(\eps^{1/2})$ range of $A$. 
The second exception is that the periodic orbit $\Gamma$ can become non-local at parameters \eJGold{close to the horizontal dashed green lines in the $(\mu,A)$ bifurcation diagrams where the slow flow has singular cycles.} In particular, period-doubling (PD) and fold of periodic orbits (LPC) bifurcations become non-local at the points marked P, and torus bifurcations (labeled NS) where a torus forms around a periodic orbit can become non-local if they are in the vicinity of the green dashed lines. \\

{\bf Global bifurcation structure.} System \eqref{resc_shnf} is highly degenerate when $B=0$. It has only one equilibrium point, which is in the fold region unless $A+C=O(\epsilon)$.  Asymptotic formulas for the first Lyapunov coefficient of the singular Hopf bifurcation when $B\neq 0$ suggest that the generalized Hopf bifurcation passes through the zero Hopf bifurcation at $B=0$ in a non-generic manner, see figures 
~\ref{fig:bdiag_Bm01C1_all_aexam_lpo_w_hom},~\ref{fig:bdiagBm01C25_all_aexam}
 and ~\ref{fig:bdiagB001C0896_all_aexam},~\ref{fig:bdiagB001C1_all_aexam},~\ref{fig:bdiagB001C16_all_aexam}
for sequences of $(\mu,A)$ diagrams in which $B$ becomes small relative to $C$. 
Also, as the generalized Hopf bifurcation passes through the zero Hopf bifurcation, the normal form type of the zero Hopf bifurcation changes~\cite{kuznetsov_book}. We conjecture that the addition of higher order terms to the system \eqref{shnf} can produce qualitative changes in its bifurcation diagram, but do not pursue that issue here.\\

We next list phenomena that produce changes in the $(\mu,A)$ bifurcation diagrams due to non-transverse intersections of these slices with bifurcation manifolds of codimensions one and two: 
\eJGold{\begin{itemize}
\item folds of generalized Hopf bifurcations
\item folds of the curves where the slow-flow has singular cycles
\item tangential intersections of the period-doubling bifurcation locus and the Hopf bifurcation loci
\item tangential intersection of the torus bifurcation locus and the Hopf bifurcation locus
\item resonant zero Hopf and generalized Hopf bifurcations
\item endpoint of the tangency curve crosses the torus bifurcation locus
\item endpoint of the tangency curve crosses the period-doubling bifurcation locus
\item the tangency curve intersects the period-doubling bifurcation locus tangentially
\item the tangency curve ``attaches'' to the Hopf curve.
\end{itemize}}

\begin{figure}[hptb]
   \begin{center}
\includegraphics[width=1.0\textwidth]{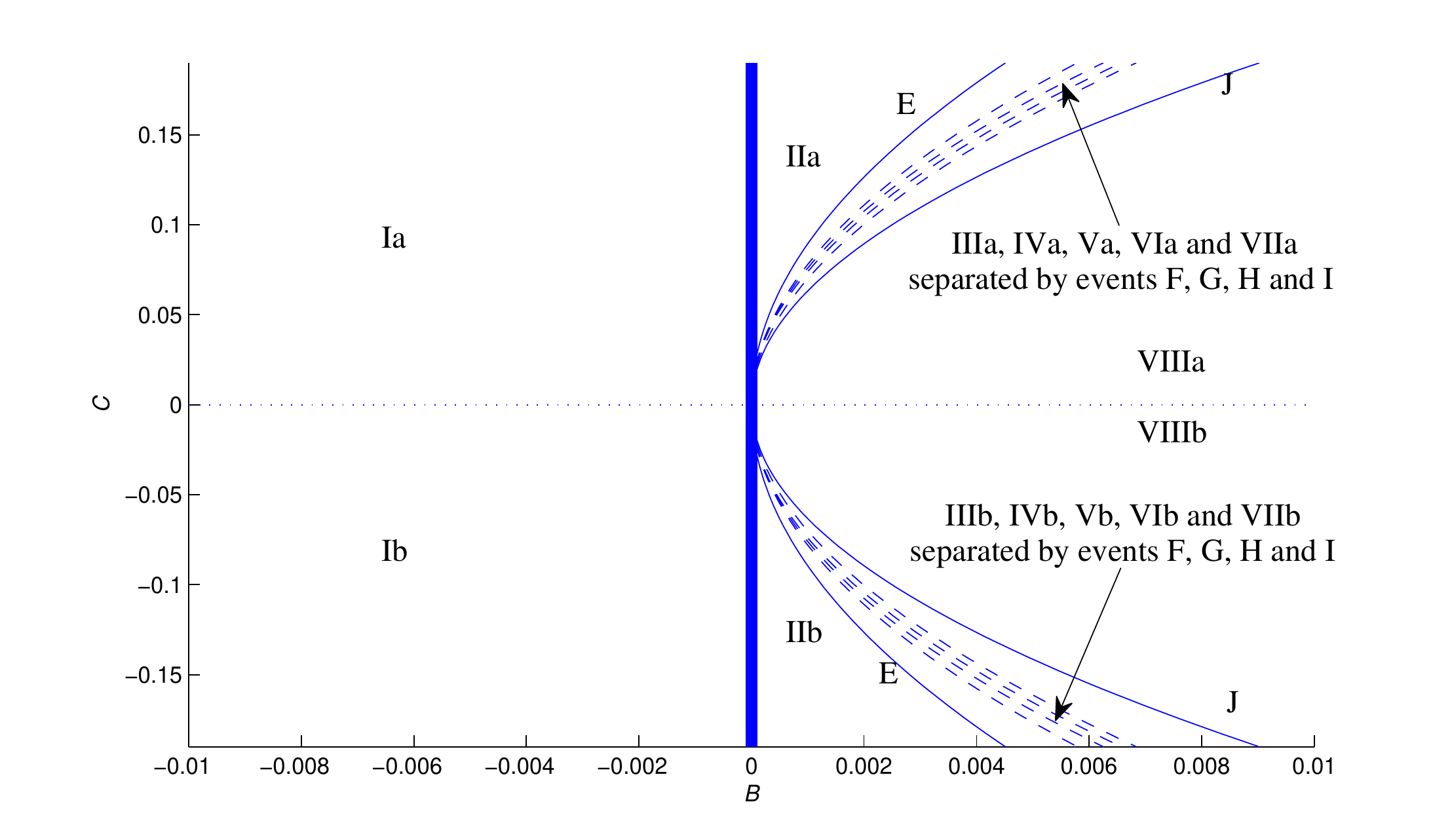} 
\caption{Regions in $(B,C)$ parameter space. Solid and dashed lines indicate the positions of codimension 3 phenomena that separate regions in $(B,C)$ parameter space yielding different corresponding $(\mu,A)$ bifurcation diagrams. }
   \label{fig:bc_regions}
\end{center}
\end{figure}

 Figure~\ref{fig:bc_regions} displays 16 regions in the $(B,C)$ parameter space \eJGold{of system ~\eqref{resc_shnf}.} These regions constitute a partitioning into regions corresponding to classes of topologically equivalent $(\mu,A)$ bifurcation diagrams, with two caveats: first, regions Ia and Ib as well as VIIIa and VIIIb have topologically equivalent $(\mu,A)$ diagrams, we merely distinguish them them for symmetry reasons. Changing the sign of parameter $C$ corresponds to reflecting bifurcation curves in $(\mu,A)$ diagrams about the $\mu$-axis, and to changing the time orientation. Under this transformation, the tangency curve transforms into a curve corresponding to a tangential intersection of $W^s(E_f)$ with $S_a$. The tangency involving $S_r$ is present in quadrants I, II and III of figure ~\ref{fig:bc_regions}, whereas the tangency involving $S_a$ is present in quadrants II,III and IV of figure ~\ref{fig:bc_regions}. Appendix B contains at least one representative $(\mu,A)$ bifurcation diagram for each of the 8 regions of $(B,C)$ parameter space in which $C>0$. 
Secondly, certain codimension three events related to the torus bifurcation and its degeneracies are not included in figure ~\ref{fig:bc_regions}. 

The degeneracy at $B=0$ (thick solid line in figure ~\ref{fig:bc_regions}) separates regions in $(B,C)$ parameter space with vastly different $(\mu,A)$ bifurcation diagrams. For $C$ fixed, $(\mu,A)$ diagrams become degenerate as $|B|$ tends to $0$ and the zero Hopf bifurcation approaches generalized Hopf bifurcations. There is only one topological type of $(\mu,A)$ diagram for \eJG{$B<0$}. The situation is more complicated for $B>0$. What follows is a description of the topological changes in the $(\mu,A)$ diagrams occurring at the curves separating the regions where $B>0$. We only treat the first quadrant of figure ~\ref{fig:bc_regions}, a description for the fourth quadrant may be obtained using the symmetry $C \mapsto -C$. \ePMnew{Note that no curves in figure ~\ref{fig:bc_regions} intersect, except maybe for $0\leq B\ll10^{-4}$. }

Regions IIa and IIIa are separated by curve E, which consists of $(B,C)$ values where the surface of generalized Hopf (GH) bifurcations is tangent to the $(\mu,A)$ slice of the parameter space. The formula for the first Lyapunov \eJGold{coefficient} indicates that the curve E is tangent to the parabola $8B=C^2$ at the origin of the $(B,C)$ plane. On this curve, there is a single generalized Hopf (GH) bifurcation of the $(\mu,A)$ bifurcation diagrams. To the left of this curve, there are two generalized Hopf bifurcations; to the right there are none. 
Curves F, G, H and I all correspond to events relating to the tangency curve or to bifurcations of periodic orbits. These curves were found by inspecting appropriate sequences of $(\mu,A)$ diagrams. A superset of the data used can be found in appendix C. 
Regions IIIa and IVa are separated by curve F, where the endpoint of the tangency curve crosses the torus bifurcation locus. In region IIIa, the tangency curve and torus bifurcation curve intersect, while they are disjoint in region IVa. 
The tangency curve coalesces with the Hopf bifurcation curve at the boundary of regions IVa and Va, represented by curve G. 
Regions Va and VIa are separated by the curve H, which consists of $(B,C)$ values where the torus bifurcation locus intersects the Hopf bifurcation locus tangentially: to the right of H, the torus bifurcation locus and the Hopf curve no longer intersect. 
In region VIa, the period-doubling bifurcation locus intersects the Hopf bifurcation in two points. As $B$ is increased, these two points coalesce when the two bifurcation loci intersect tangentially on the curve I, after which they no longer meet in regions VIIa and VIIIa. 
Regions VIIa and VIIIa are separated by the curve J where the two $(\mu,A)$ curves of homoclinic bifurcations coalesce as $B$ increases. \eJGold{The curve J is approximated by the quadratic equation $C^2-4B=0$.}\\

{\bf Torus bifurcation and the global bifurcation structure. } For $B<0$, taking account of the torus bifurcation curve and its degeneracies (resonances) does not yield more topologically different $(\mu,A)$ bifurcation diagrams. For $B>0$, the division of $(B,C)$ parameter space into regions with topologically equivalent $(\mu,A)$ diagrams gets significantly more complicated and is not presented here. The paragraphs of section 4 headed ``Periodic orbit bifurcations'' give a general description of the positions of torus bifurcation curves relative to other bifurcation curves. Data on the position of the crossing of a 1:2 resonance with the Hopf curve can be found in appendix C. \\

{\bf One dimensional bifurcation diagrams.}
Table  ~\ref{table: sequence_table} lists the sequences of local bifurcations shown in the $(\mu,A)$ diagrams of this paper that occur as the primary bifurcation parameter $\mu$ is increased. 
We list only sequences occurring well within the $(B,C)$ regions we have identified. There could be some additional sequences, e.g. \ePM{extremely close to the zero Hopf bifurcation, or} where $|B|<<C^2$, i.e. where the zero Hopf, generalized Hopf and P bifurcations occur at very similar values of ~$A$. 
Two branches of periodic orbits emanate from LPC bifurcation points. \ePM{If a sequence begins with an LPC bifurcation and also contains a torus bifurcation or period-doubling, then one branch of periodic orbits emanating from the LPC bifurcation has a torus bifurcation and period-doubling, while the other one ends in a singular Hopf bifurcation}. 
$W^u(E_f)$ and $S_r$ intersect after the tangency in the bifurcation sequences. 

\begin{table}

\begin{center} \footnotesize

\begin{tabular}{ l | l | l }
Sequence \#	&Sequence description	& $(B,C)$ \ePM{region} present \\
\hline
1	&H\textsubscript{sup}						&Ia, IIa, IIIa, IVa, Va, Ib, IIb, IIIb, IVb, Vb\\
2	&H\textsubscript{sup} - (SN)					&IIa, IIIa, IVa, Va, IIb, IIIb, IVb,Vb\\
3	&(SN) - H\textsubscript{sup} - T $\pm$ PD		&Ib\\
4	&(SN) - H\textsubscript{sup} - T $\pm$ NS - PD		&Ia\\
5	&H\textsubscript{sup} - LPC					&Ia, Ib\\
6	&H\textsubscript{sup} - T $\pm$ PD				&Ia, IIa, Ib\\
7	&H\textsubscript{sup} - T $\pm$ NS - PD			&IIa\\
8	&LPC - H\textsubscript{sup} - NS - PD			&IIb, IIIb\\
9	&LPC - H\textsubscript{sup} - PD				&IIb, IIIb\\
10	&LPC - H\textsubscript{sup} - T $\pm$ PD			&IIa\\
11	&LPC - NS - H\textsubscript{sup}  - PD			&IIb, IIIb\\
12	&LPC - PD - H\textsubscript{sup}				&IIb, IIIb, IVb\\
13	&H\textsubscript{sub} 						&Ia, IIa, IIIa, IVa, Va, Ib, IIb, IIIb, IVb, Vb\\
14	&H\textsubscript{sub} - (SN)					&IIa, IIIa, IVa, Va, IIb, IIIb, IVb, Vb\\
15	&(SN) - H\textsubscript{sub} - PD				&Ia\\
16	&(SN) - H\textsubscript{sub} - NS - PD			&Ib\\
17	&H\textsubscript{sub} - LPC					&Ia, Ib\\
18	&H\textsubscript{sub} - PD					&Ia, Ib, IIb\\
19	&H\textsubscript{sub} - NS - PD				&IIb\\
20	&LPC - H\textsubscript{sub} - NS - PD			&IIa, IIIa\\
21	&LPC - H\textsubscript{sub} - T $\pm$ NS - PD		&IIa, IIIa  \\
22	&LPC - H\textsubscript{sub} - T $\pm$ PD 		&IIa, IIIa \\
23	&LPC - H\textsubscript{sub}  - PD				&IIa, IIIa\\
24	&LPC - NS - H\textsubscript{sub}  - T $\pm$ PD	&IIa, IIIa\\
25	&LPC - PD - H\textsubscript{sub} 				&IIa, IIIa, IVa
\end{tabular}
\caption{List of bifurcation sequences occurring in system \eqref{resc_shnf}  as the main bifurcation parameter $\mu$ is increased. 
Each sequence is described in terms of the codimension 1 bifurcations involved, $\mu$ increasing. 
The sign ``-'' separates codimension 1 bifurcations occurring at different values of $\mu$. 
\ePM{A ``$\pm$'' sign appears between two bifurcations where the first bifurcation listed either occurs before the second or indistinguishably close to the second, at the resolution of the diagrams in this paper. 
Saddle-node bifurcations are listed in parentheses, since they are only local close to the zero Hopf point. }}
\label{table: sequence_table}
\end{center}
\end{table}

\section{An Example: the Koper Model}

\begin{figure}[hptb]
\begin{center}
\includegraphics[width=1.0\textwidth]{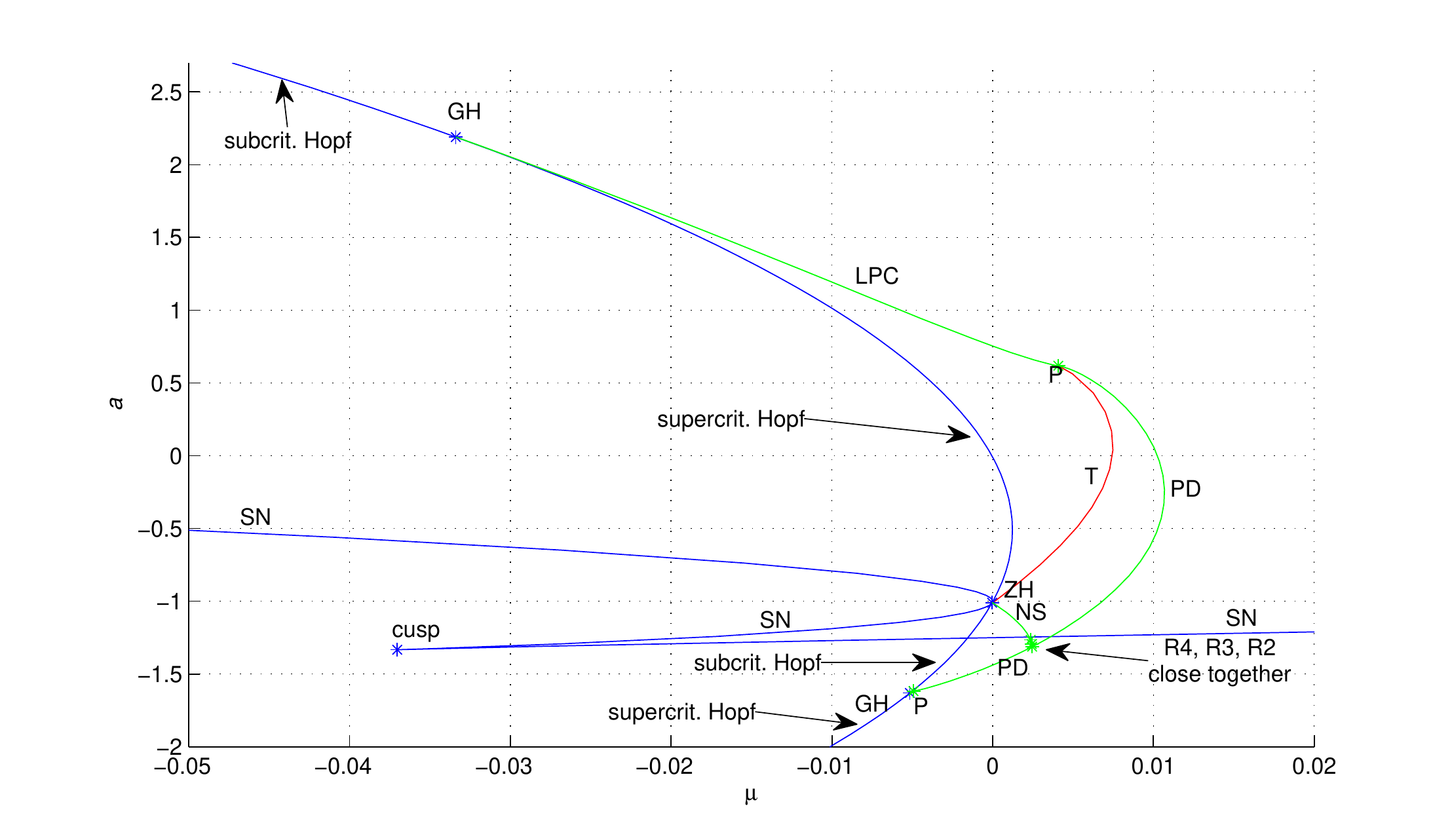} 
\caption{Bifurcation diagram for cubic singular Hopf normal form, $B=-0.01,C=0.1$.}
\label{fig:NEW_bdiag_eps01Bm1C1_param01}
\end{center}
\end{figure}

This section demonstrates how our normal form analysis can be used in investigations of systems undergoing singular Hopf bifurcation. We use a vector field that was studied first by Koper~\cite{koper_paper} as a model of mixed mode oscillations in chemical systems:
\ePMnew{\begin{equation} 
\begin{split}
\eps_1 \,\dot{x} & = 
k\,y-x^3 +3\,x-\lambda\\
\dot{y} & = 
x-2\,y+z \\
\dot{z} & = 
\eps_2\,(y-z)\\
\end{split}
\label{koper}
\end{equation}}
Further discussion of this model can be found in \ePMnew{~\cite{KrupaPopovicKopell}, \ePMnew{~\cite{Kuehn}} as well as} Desroches  et al.~\cite{mmo_paper} and references cited there. Using affine coordinate changes and a time rescaling (see ~\cite{mmo_paper}), the Koper vector field can be written in the following form: 
\begin{equation} 
\begin{split}
\dot{x}
& = 
(y-x^2-x^3)/\eps  \\
\dot{y} & = 
z-x \\
\dot{z}& =
-\mu-a \,x -b \,y -c\, z\\
\end{split}
\label{cubic_shnf}
\end{equation}
Note that 
\begin{itemize}
 \item 
the system \eqref{cubic_shnf} has one more parameter than the system \eqref{koper}, so the Koper vector field is a subfamily of \eqref{cubic_shnf}, and
\item
the system \eqref{cubic_shnf} is a variant of the singular Hopf normal form \eqref{shnf} in which the parabolic critical manifold has been replaced by a cubic critical manifold ~\cite{mmo_paper} by the addition of a single, higher order term to the equation.
\end{itemize}
With the scaling of parameters and coordinates used to obtain system~\eqref{resc_shnf} from system~\eqref{shnf},  \eqref{cubic_shnf} transforms to 
\begin{equation} 
\begin{split}
\dot{X}
& = 
(Y-X^2-\eps^{1/2}X^3) \\
\dot{Y} & = 
Z-X \\
\dot{Z}& =
-\mu-A X -B Y -C Z\\
\end{split}
\label{resc_cubic_shnf}
\end{equation}
It is apparent that the system~\eqref{resc_cubic_shnf} is a small perturbation of the system~\eqref{cubic_shnf} in a bounded region of $(X,Y,Z)$ space.
Figures~\ref{fig:NEW_bdiag_eps01Bm1C1_param01} and \ref{fig:bdiag_Bm01C1_all_aexam_lpo_w_hom} show $(\mu,A)$ bifurcation diagrams for systems~\eqref{resc_cubic_shnf} and \eqref{resc_shnf} at the same values of $(B,C) = (-0.01,0.1)$. The similarity of these two figures supports the use of the normal form \eqref{shnf} to study local bifurcations near a singular Hopf bifurcation. Note that the cusp of the saddle-node curve in figure~\ref{fig:NEW_bdiag_eps01Bm1C1_param01} is not local: it tends to   infinity as $\eps \to 0$.

\begin{figure}[hptb]
\begin{center}
\includegraphics[width=1.0\textwidth]{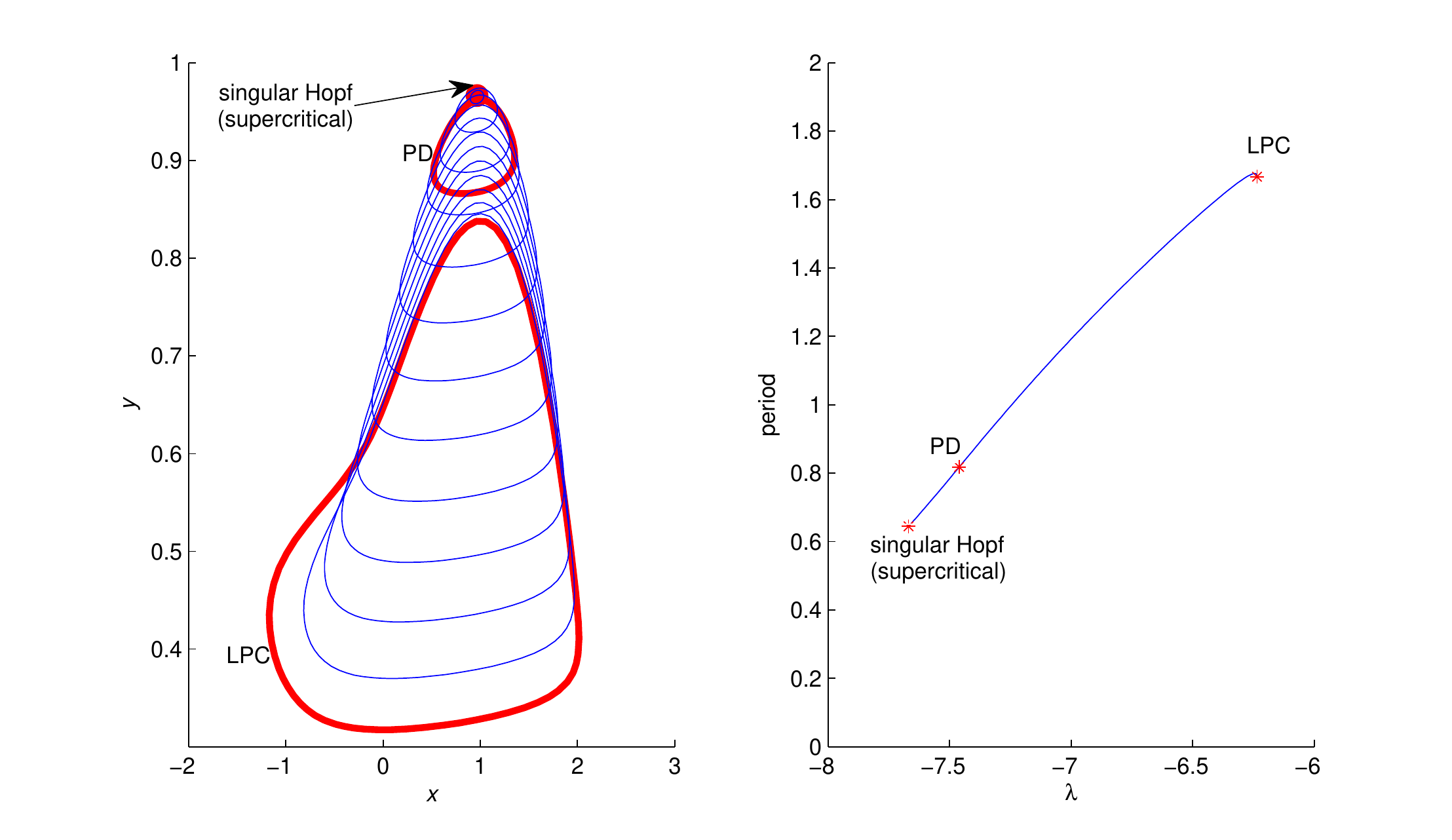} 
\caption{Periodic orbits emanating from a supercritical Hopf bifurcation for system ~\eqref{koper} and  $\eps_1=0.1,k=-10$. Observe that the periodic orbits follow the cubic shape of the (parameter-dependent) critical manifold more as $\lambda$ increases.}
\label{fig:plotting_koper_periodic_orbits}
\end{center}
\end{figure}

Consider the Koper model with $\eps_2=1$. When additionally $\eps_1=0.1$ and $k=-10$, there is a supercritical singular Hopf bifurcation at approximately  $\lambda \approx-7.670$. Continuing this family of periodic orbits with varying lambda using AUTO, we find that the emerging periodic orbit undergoes a period-doubling bifurcation at $\lambda \approx-7.461$, and the ``undoubled'' periodic orbit undergoes a fold of periodic orbits (LPC) bifurcation at $\lambda \approx-6.235$. The periods of the periodic orbits at the bifurcation points are approximately $0.64$, $0.82$ and $1.67$ respectively. If $\tau$ denotes the initial period of the periodic orbits emerging from the singular Hopf bifurcation, local bifurcations of periodic orbits typically occur with periods less than $1.5\tau$. This suggests that the LPC bifurcation is not local. Figure ~\ref{fig:plotting_koper_periodic_orbits} shows that indeed, the shape of the periodic orbit at the LPC bifurcation follows the cubic shape of the critical manifold, rather than staying close to the fold curve at which the singular Hopf bifurcation occurred. \\\

\ePM{Among the sequences with a supercritical Hopf bifurcation, only 3, 4 and 6 to 12 also have a period-doubling. Note that sequences 8 to 12 have a period-doubling bifurcation and a supercritical Hopf bifurcation on distinct branches of periodic orbits, and that sequences 4 and 7 contain torus bifurcations. Thus, only sequences 3 and 6 match the continuation data. This suggests that there is a tangency of invariant manifolds either between the supercritical Hopf bifurcation and the period-doubling bifurcation or very close to the period-doubling bifurcation:} indeed, numerical calculations show that \ePMnew{the tangency} occurs at $\lambda \approx-7.539$. 

\begin{figure}[hptb]
\begin{center}
\includegraphics[width=1.0\textwidth]{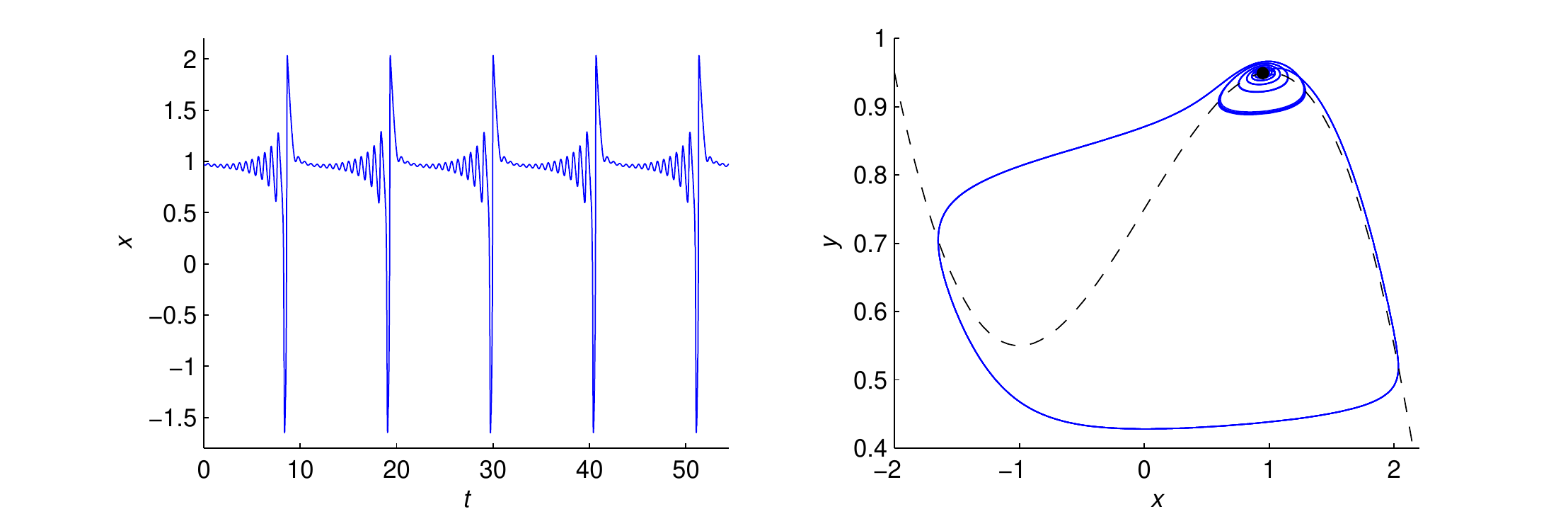} 
\caption{Mixed mode oscillations in system ~\eqref{koper} with  $\eps_1=0.1,\eps_2 = 1, k=-10, \lambda = -7.50$: the left panel shows the time-series of the $x$-coordinate of an MMO trajectory, the right panel shows the same trajectory in the $xy$-plane. $E_f$ is marked with a black dot, the critical manifold is drawn with a dashed black line. }
\label{fig:koper_mmo}
\end{center}
\end{figure}

Due to the S-shape of the critical manifold in the Koper model, trajectories in the unstable manifold in the vicinity of the fold region that escape from the fold region can return to the fold region. As a result, the tangency of invariant manifolds can mark the onset of mixed-mode oscillations:  a trajectory approaches the equilibrium point near the fold region, then follows the equilibrium's two-dimensional unstable manifold, leaves the fold region, and approaches the equilibrium once more along its stable manifold. The Koper system has mixed mode oscillations at $\lambda=-7.50$, illustrated in figure ~\ref{fig:koper_mmo}. Thus local analysis of system \eqref{resc_shnf} led to information about the Koper model as an alternative to applying numerical continuation tools directly on the Koper model. This approach can be used to explore the effects of varying additional parameters in the Koper model as well as to study other systems in which singular Hopf bifurcation occurs. \\

\section{Discussion}

This paper includes a comprehensive analysis of of bifurcations that occur in the normal form~\eqref{shnf} for a singular Hopf bifurcation. In particular, we analyze bifurcations of equilibrium points, periodic orbits originating from these equilibria and ``simple'' tangencies between the unstable manifold \ePMnew{$W^u(E_f)$} of an equilibrium \ePMnew{$E_f$} and a repelling slow manifold \ePMnew{$S_r$} of the system. We approach the study of these bifurcations in a hierarchical fashion, designating the parameter $\mu$ as the primary bifurcation parameter and then performing analyses of parallel two dimensional slices throughout the four dimensional parameter space of the system. While our analysis is comprehensive, it is hardly complete. We point to four aspects of the bifurcations of this system for further study.

In some parameter ranges, the system~\eqref{shnf} has complicated invariant sets in $O(\sqrt{\eps})$ neighborhoods of the origin. Torus and period-doubling bifurcations of periodic orbits give rise to invariant sets that undergo further bifurcations. We have not pursued the study of these invariant sets or their bifurcations. Trajectories in these sets typically make multiple oscillations when projected to the $(x,y)$ coordinate plane. Similarly, there are more complex intersections of invariant manifolds than those we have investigated in this paper. The three two dimensional invariant manifolds that we have studied are the attracting and repelling slow manifolds and the unstable manifold of an equilibrium point. In some parameter regimes, the attracting and repelling slow manifolds have multiple intersections that constitute maximal canards of a folded node point of the system~\cite{B,GHai,Wech}. These manifolds and their intersections bound rotational sectors of trajectories that make different numbers of oscillations as they pass through the folded node region. The unstable manifold of the equilibrium can also become entangled with the repelling slow manifold and have multiple intersections that separate regions of the unstable manifold that flow to $x=-\infty$ after making different numbers of oscillations. In systems like the Koper model where trajectories jumping from the vicinity of a singular Hopf bifurcation make a global return to that region, tangencies between the different invariant manifolds play a significant role in the bifurcations of MMOs. We have not considered here secondary tangencies of the unstable manifold and the repelling slow manifold that are analogous to the bifurcations of secondary canards of a folded node analyzed by Wechselberger~\cite{Wech}.

The second aspect of singular Hopf bifurcation that calls for further investigation is the presence of highly degenerate bifurcations within the normal form. In particular, we have seen that the system is very degenerate when $b=0$. In analogy with the study of codimension two bifurcations of equilibria~\cite{GH}, the addition of higher order terms could resolve this degeneracy. We do not study that possibility here. Even in systems with a single time scale, the codimension three bifurcations that correspond to degenerate zero-Hopf bifurcations have not been studied thoroughly. 

\begin{figure}[hptb]
\begin{center}
\includegraphics[width=1.0\textwidth]{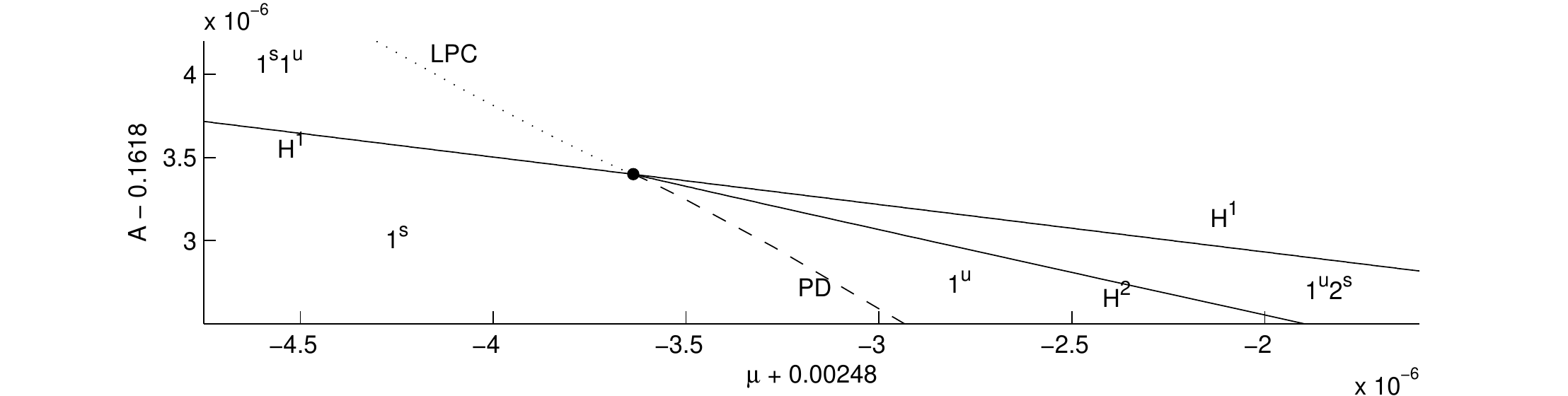} 
\caption{Bifurcation diagram showing bifurcation curves of periodic orbits in the vicinity of homoclinic orbit flip at $(\mu,A,B,C)=(-0.00248364, 0.161803, -0.01, -0.1)$. Numbers \eJGold{with superscripts $s$ or $u$} indicate the number of loops as well as the stability of periodic orbits in each region of parameter space. Unstable periodic orbits are of saddle-type. The labels $\mathrm{H^1}$ and $\mathrm{H^2}$ distinguish homoclinic bifurcation loci of periodic orbits with one and two loops respectively.}
   \label{plot_LPC_PD_H2_H1_hom_orbit_flip_caseB_strange}
\end{center}
\end{figure}

Third, our investigation of the the points $P$ where period-doubling bifurcation, fold of periodic orbits bifurcation and the slow flow has singular cycles meet at the degenerate homoclinic bifurcation remains incomplete. These points appear to be homoclinic orbit flips or inclination flips~\cite{n_cont_dyn_sys}, but we have not reconciled aspects of the behavior found here with that predicted by the unfoldings of generic orbit flips or inclination flips. Investigations of phase portraits suggest that if $B<0$, there is a homoclinic orbit flip at approximately $(\mu,A)=(B/4,(-C \pm \sqrt{C^2-4B})/2)$, provided that the square root evaluates to a real number. However, we find that the relative positions of the bifurcation curves are interchanged from those predicted by the theory, cf. \cite{death_of_period_doublings} and figure ~\ref{plot_LPC_PD_H2_H1_hom_orbit_flip_caseB_strange}.
If $B>0$, we surmise the existence of an  homoclinic orbit [inclination] flip of type C~\cite{death_of_period_doublings} at approximately $(\mu,A)=(B/4,(-C - \sqrt{C^2-4B})/2)$ if the parameter $C$ is positive [negative], and at $(\mu,A)=(B/4,(-C + \sqrt{C^2-4B})/2)$ if parameter $C$ is negative [positive], cf. section 3.  \\

\newpage
\appendix

\section{Bifurcation labels}

\begin{table}[h]

\begin{center} \footnotesize

\begin{tabular}{ l | l | l}
Bifurcation label 	& Bifurcation / explanation	& Codimension \\
\hline
H\textsubscript{sup}	& supercritical Hopf bifurcation		& 1	\\
H\textsubscript{sub}	&subcritical Hopf bifurcation	& 1	\\
SN		&saddle-node bifurcation 			& 1	\\
PD		&period-doubling bifurcation		& 1	\\
NS		&torus bifurcation  / Neimark-Sacker bifurcation & 1 \\
LPC		&fold of periodic orbits bifurcation  / limit point of cycles bifurcation  & 1\\
H		& homoclinic bifurcation 	& 1	\\
T		& tangency of invariant manifolds bifurcation	& 1	\\
S		& approximate symmetry of the slow flow  & 1 \\
ZH		& zero Hopf bifurcation  / fold Hopf bifurcation 	& 2 \\
GH		& generalized Hopf bifurcation  / Bautin bifurcation / Gavrilov-Guckenheimer bifurcation & 2 \\
P		& see sections 4 and 5	& 2 \\
R1		& 1:1 resonance on a torus bifurcation curve 	& 2 	\\
R2		& 1:2 resonance on a torus bifurcation curve	& 2 	\\
R3		& 1:1 resonance on a torus bifurcation curve 	& 2 \\
R4		& 1:4 resonance on a torus bifurcation curve	& 2 \\
E		& fold of generalized Hopf bifurcations	& 3 \\
F		& endpoint of the tangency curve crosses the torus bifurcation curve	& 3 \\
G		& tangency curve coalesces with the Hopf bifurcation curve	& 3 \\
H		& torus bifurcation curve intersect the Hopf bifurcation curve tangentially	& 3 \\
I		& period-doubling bifurcation curve intersects the Hopf bifurcation curve tangentially	& 3 \\
J		& fold of the curves where the slow-flow has singular cycles	& 3
\end{tabular}
\caption{Table showing the abbreviations for bifurcations used in this paper, together with the full name or description of the bifurcation, as well as their codimension. Different names for the same bifurcation are separated by slashes. }
\label{table:labelexplanation}
\end{center}
\end{table}

\newpage

\section{Two dimensional bifurcation diagrams}

This appendix is a catalog of $(\mu,A)$ bifurcation diagrams for system ~\eqref{resc_shnf}, including at least one sample diagram for each of the $(B,C)$ regions Ia to VIIIa. 


\begin{figure}[hptb]
\begin{center}
\includegraphics[width=1.0\textwidth]{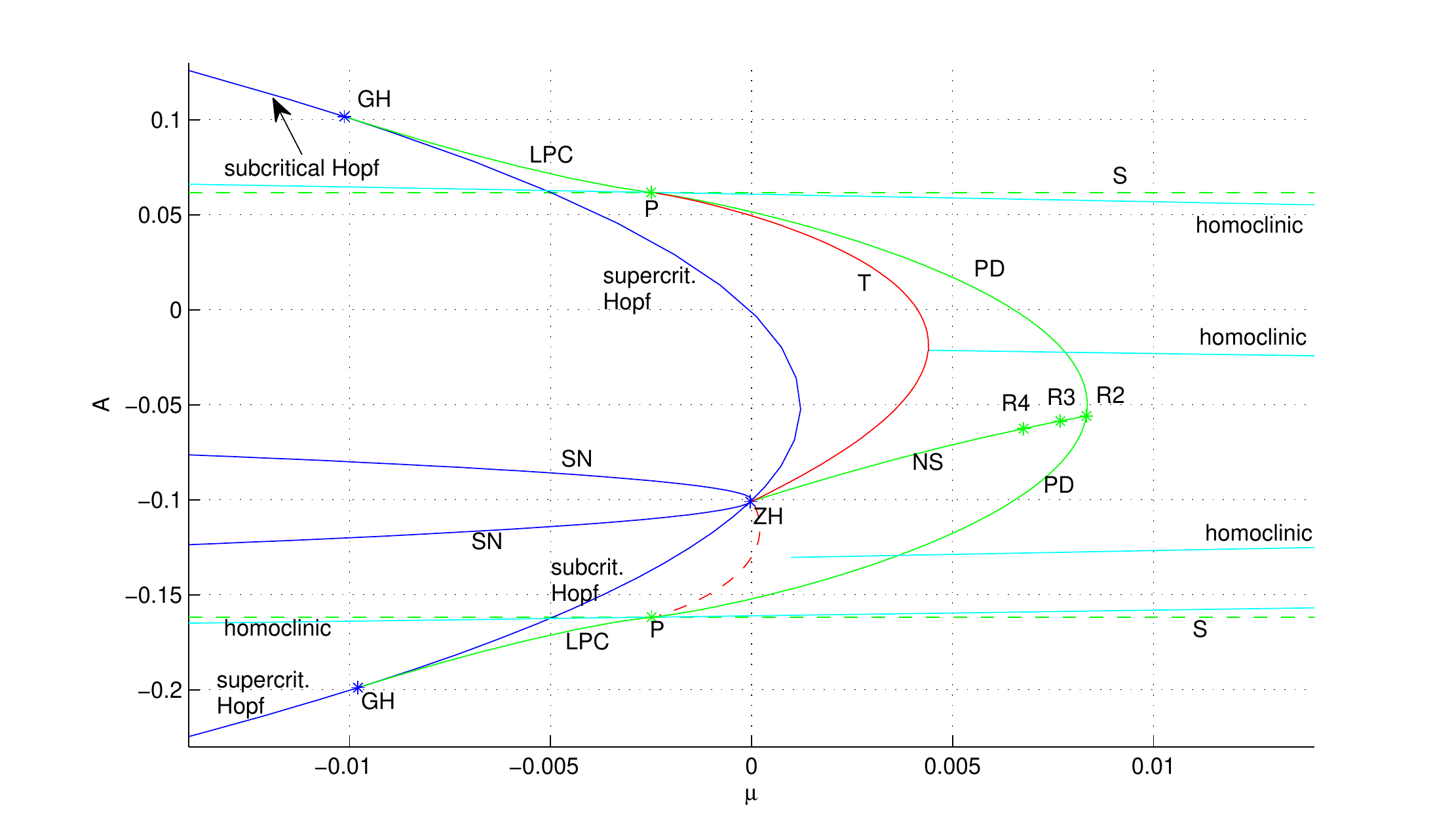} 
\caption{Region Ia: $(\mu,A)$ bifurcation diagram for $(B,C)=(-0.01,0.1)$. 
This diagram shows homoclinic bifurcations as cyan curves. The top and bottom homoclinic curves have equilibria close the origin, while the middle two have equilibria distant from the origin. The first and the third homoclinic curves connect for $O(\epsilon^{1/2})$ values of $\mu$, as do the second and fourth homoclinic curves. 
We conjecture that the endpoint of the top curve is close to where the periodic orbit born in the singular Hopf bifurcation ceases to be the only $\omega$ limit set for trajectories in the two-dimensional unstable manifold of the singular Hopf equilibrium, and that the endpoint of the bottom curve is close to where the periodic orbit born in the singular Hopf bifurcation ceases to be the only $\alpha$ limit set for trajectories in the two-dimensional stable manifold of the singular Hopf equilibrium. The location of a tangential intersection of $S_a$ with $W^s(E_f)$ is drawn with a dashed red line. 
}
\label{fig:bdiag_Bm01C1_all_aexam_lpo_w_hom}
\end{center}
\end{figure}

\begin{figure}[hptb]
\begin{center}
\includegraphics[width=1.0\textwidth]{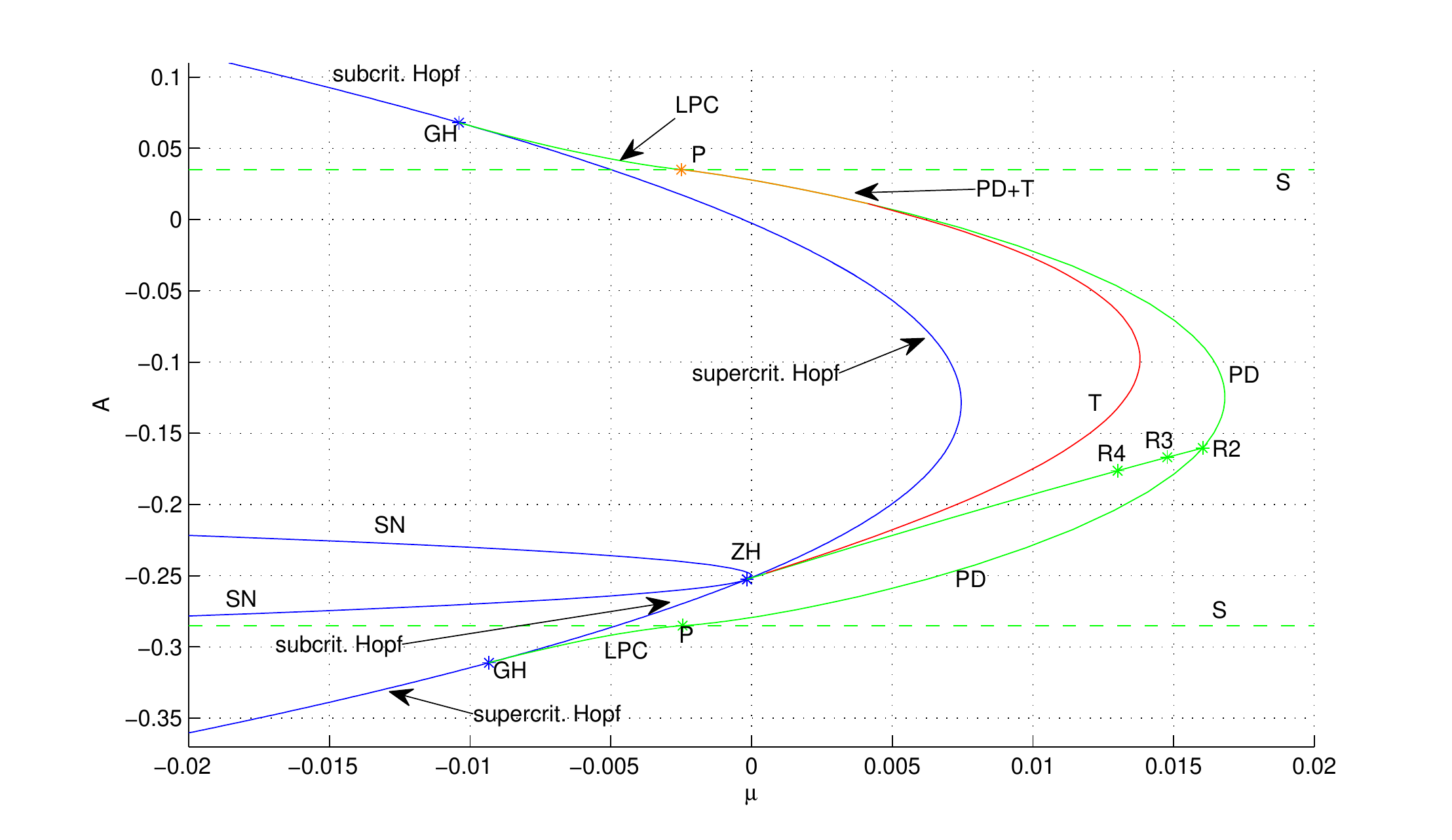} 
\caption{Region Ia: $(\mu,A)$ bifurcation diagram for $(B,C)=(-0.01,0.25)$. Observe that the tangency curve in some regions nearly coincides with the period-doubling and the torus bifurcation curves. This is not the case for all $(B,C)$ in region Ia. For $B<0$ with large absolute value or $C$ less positive, the torus bifurcation curve coincides with with the tangency curve on shorter segments, or not at all. For $B<0$ with small absolute value, the tangency curve tends to separate from the period-doubling curve closer to the point labeled P. }
\label{fig:bdiagBm01C25_all_aexam}
\end{center}
\end{figure}

\begin{figure}[hptb]
\begin{center}
\includegraphics[width=1.0\textwidth]{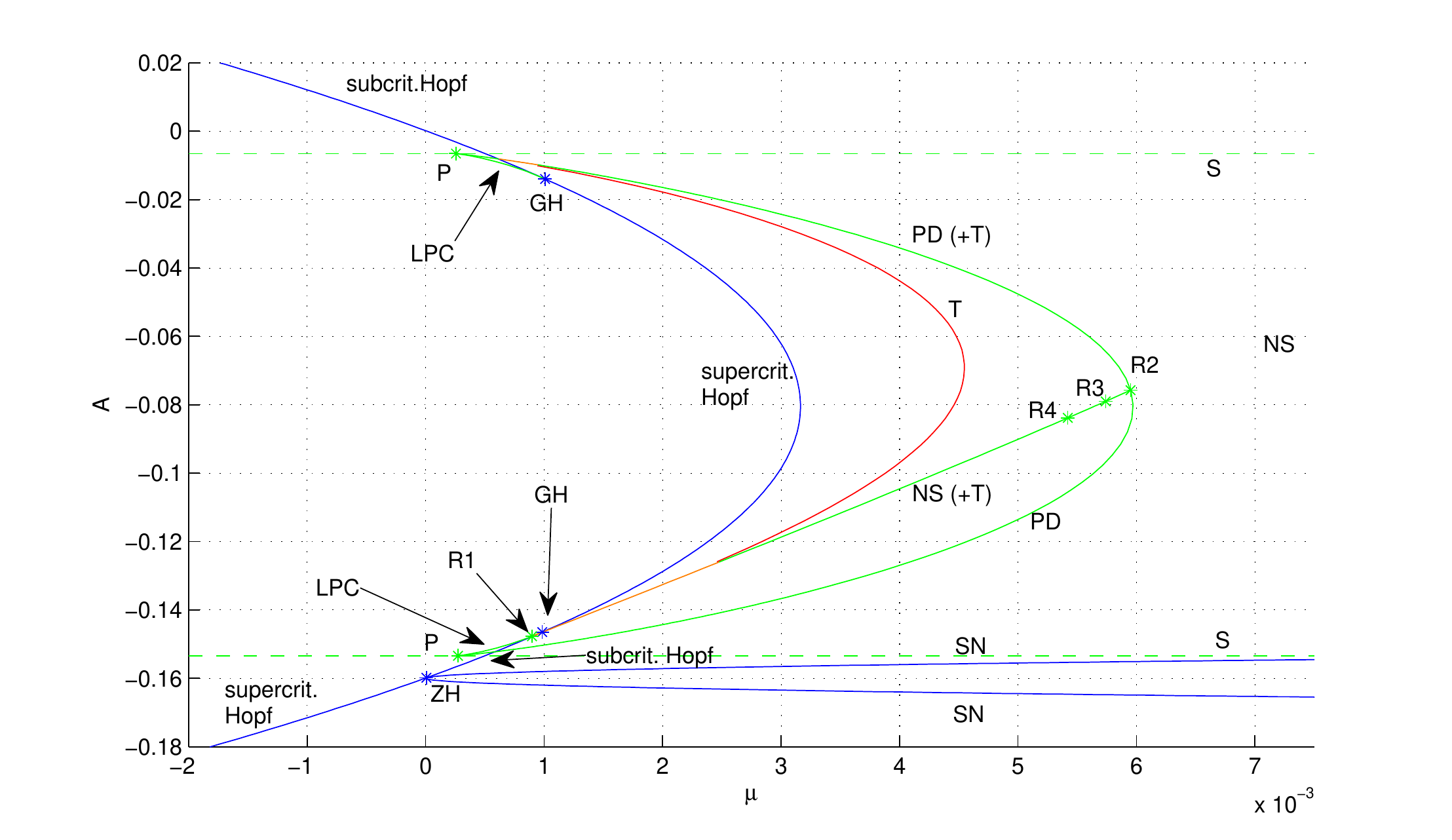} 
\caption{Region IIa: $(\mu,A)$ bifurcation diagram for $(B,C)=(0.001,0.16)$. Note that the tangency curve partly coincides with the period-doubling curve, cf. figure ~\ref{fig:bdiagB001C1_all_aexam}. This is typically observed for very large values of $C$, cf. figure ~\ref{fig:codim3_events}. }
\label{fig:bdiagB001C16_all_aexam}
\end{center}
\end{figure}

\begin{figure}[hptb]
\begin{center}
\includegraphics[width=1.0\textwidth]{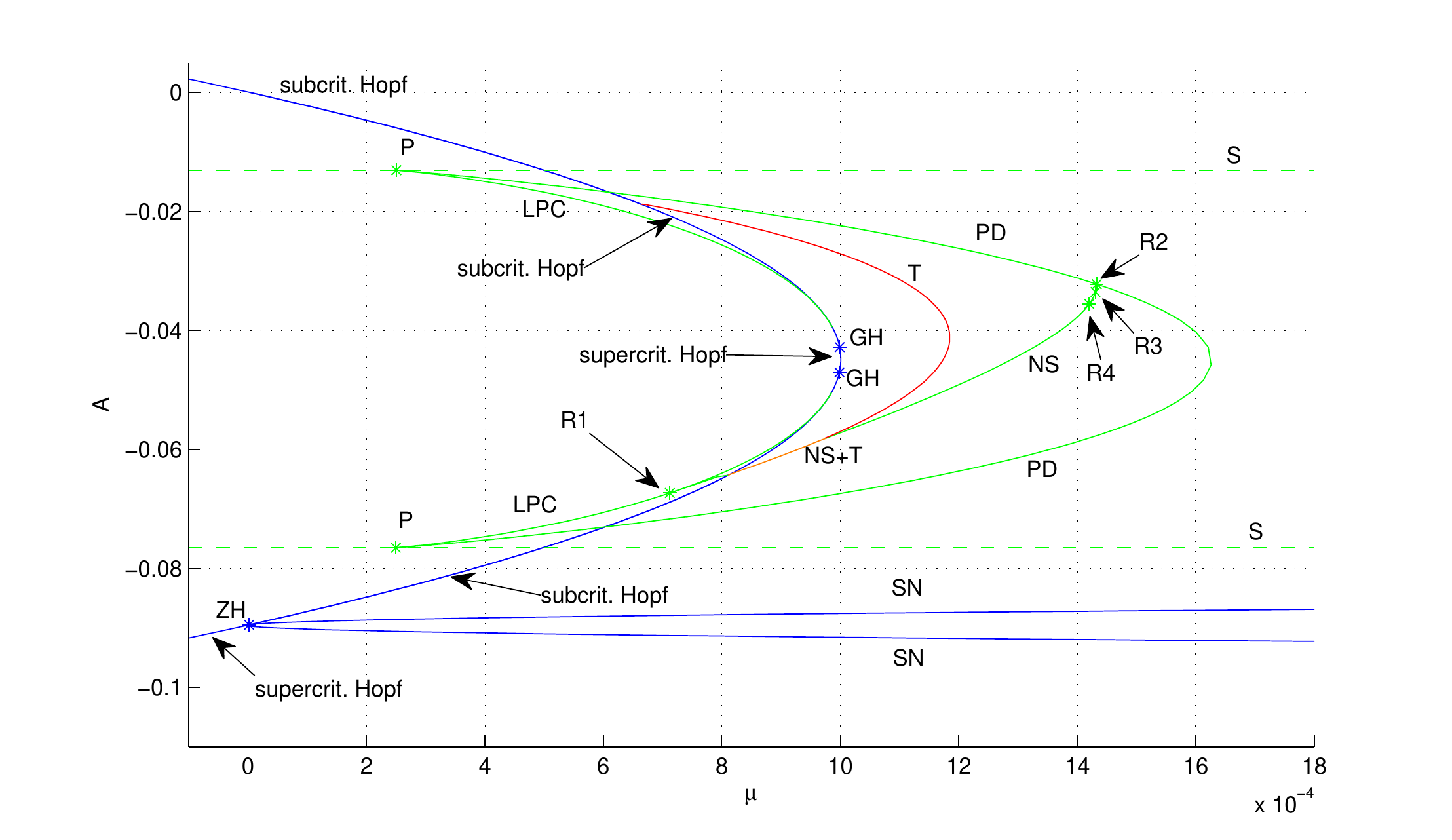} 
\caption{Region IIa: $(\mu,A)$ bifurcation diagram for $(B,C)=(0.001,0.0896)$. The generalized Hopf bifurcations are close together, and vanish at $(B,C)=(0.001, 0.0894)$.}
\label{fig:bdiagB001C0896_all_aexam}
\end{center}
\end{figure}

\begin{figure}[hptb]
\begin{center}
\includegraphics[width=1.0\textwidth]{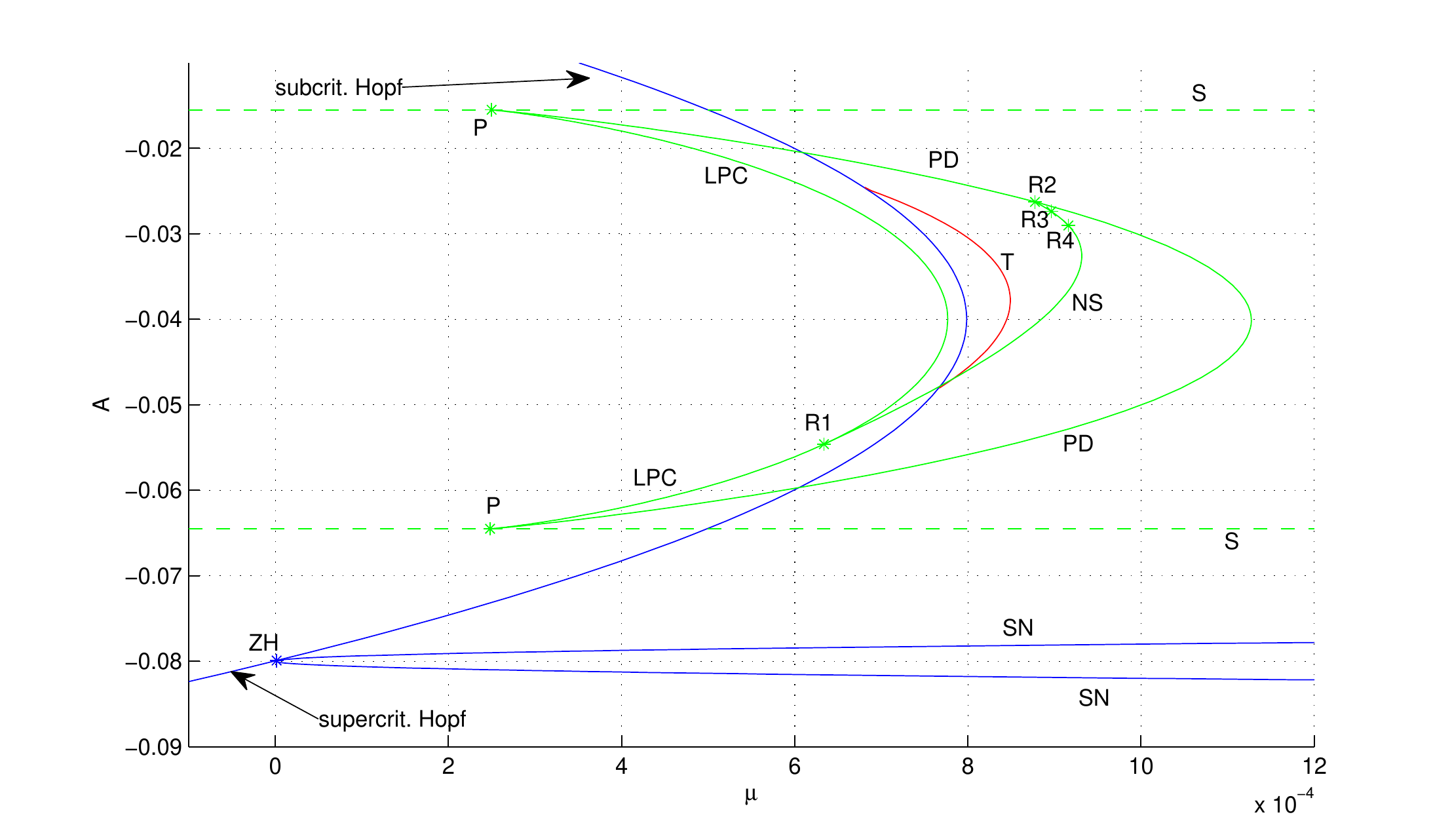} 
\caption{Region IIIa: $(\mu,A)$ bifurcation diagram for $(B,C)=(0.001,0.08)$.}
\label{fig:bdiagB001C08_all_aexam}
\end{center}
\end{figure}

\begin{figure}[hptb]
\begin{center}
\includegraphics[width=1.0\textwidth]{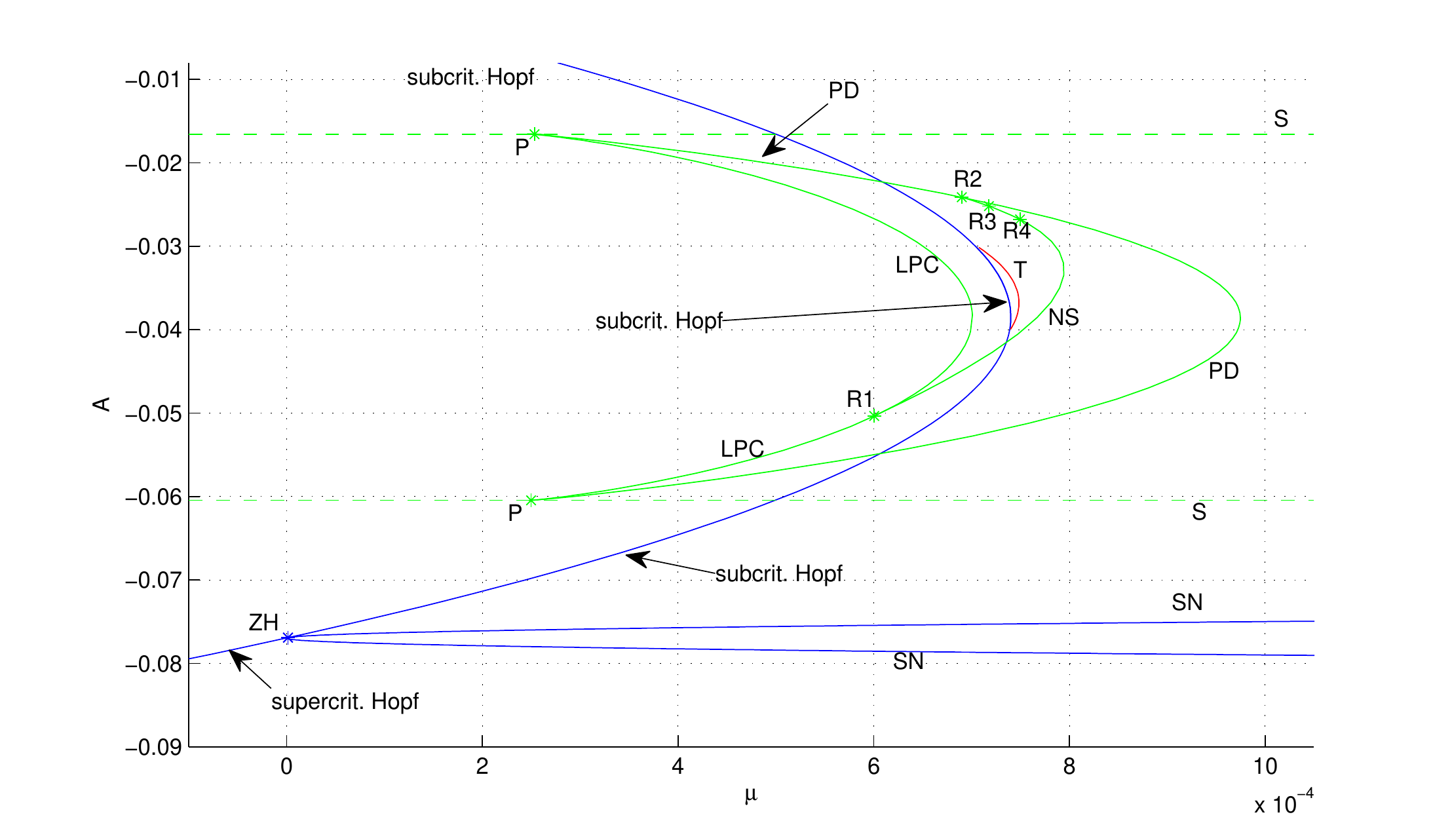} 
\caption{Region IVa: $(\mu,A)$ bifurcation diagram for $(B,C)=(0.001,0.077)$. The tangency curve and the torus bifurcation curve do not meet in this diagram.}
\label{fig:bdiagB001C077_all_aexam}
\end{center}
\end{figure}

\begin{figure}[hptb]
\begin{center}
\includegraphics[width=1.0\textwidth]{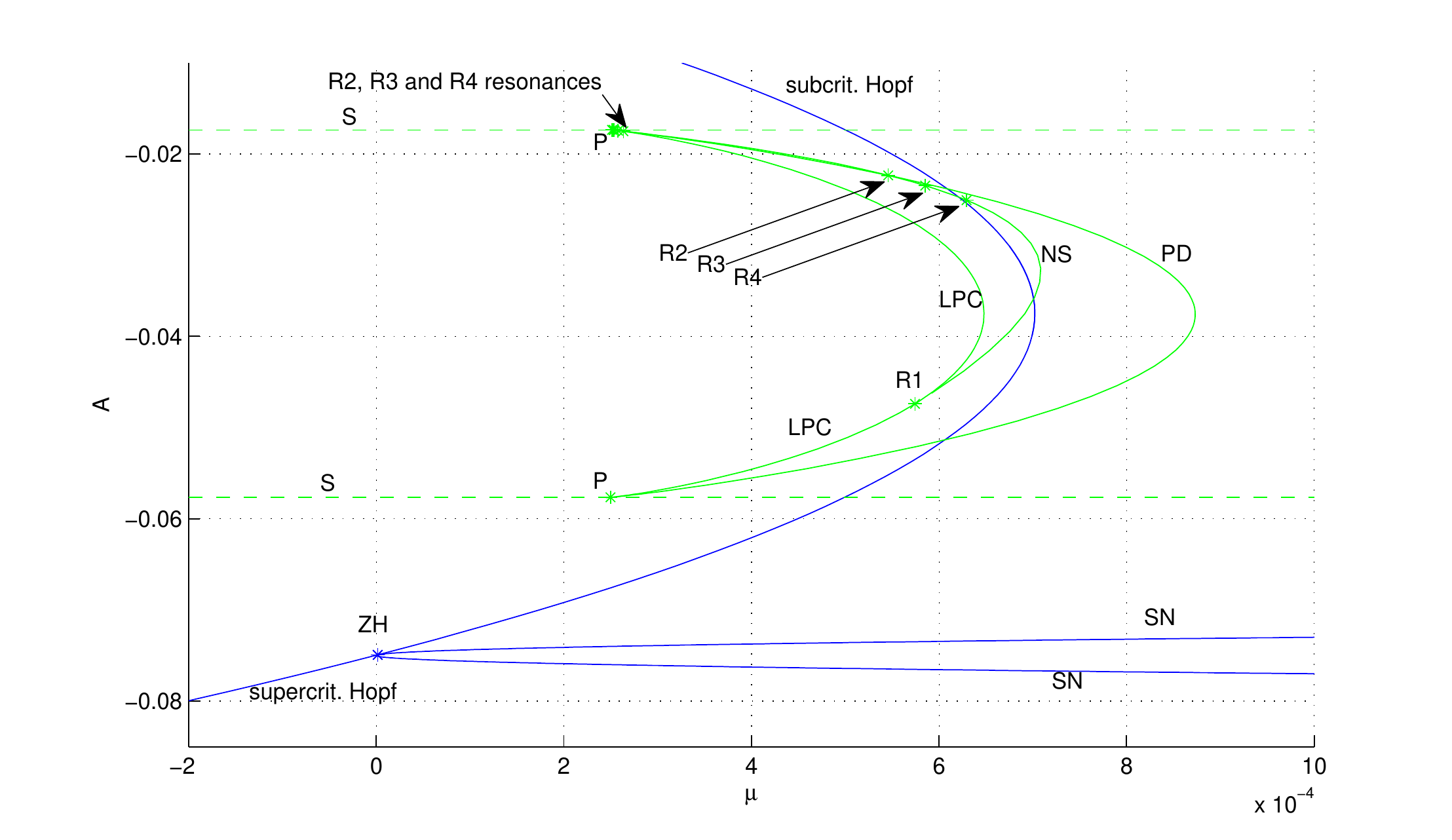} 
\caption{Region Va: $(\mu,A)$ bifurcation diagram for $(B,C)=(0.001,0.075)$.}
\label{fig:NEW_bdiag_B001C075_param28}
\end{center}
\end{figure}

\begin{figure}[hptb]
\begin{center}
\includegraphics[width=1.0\textwidth]{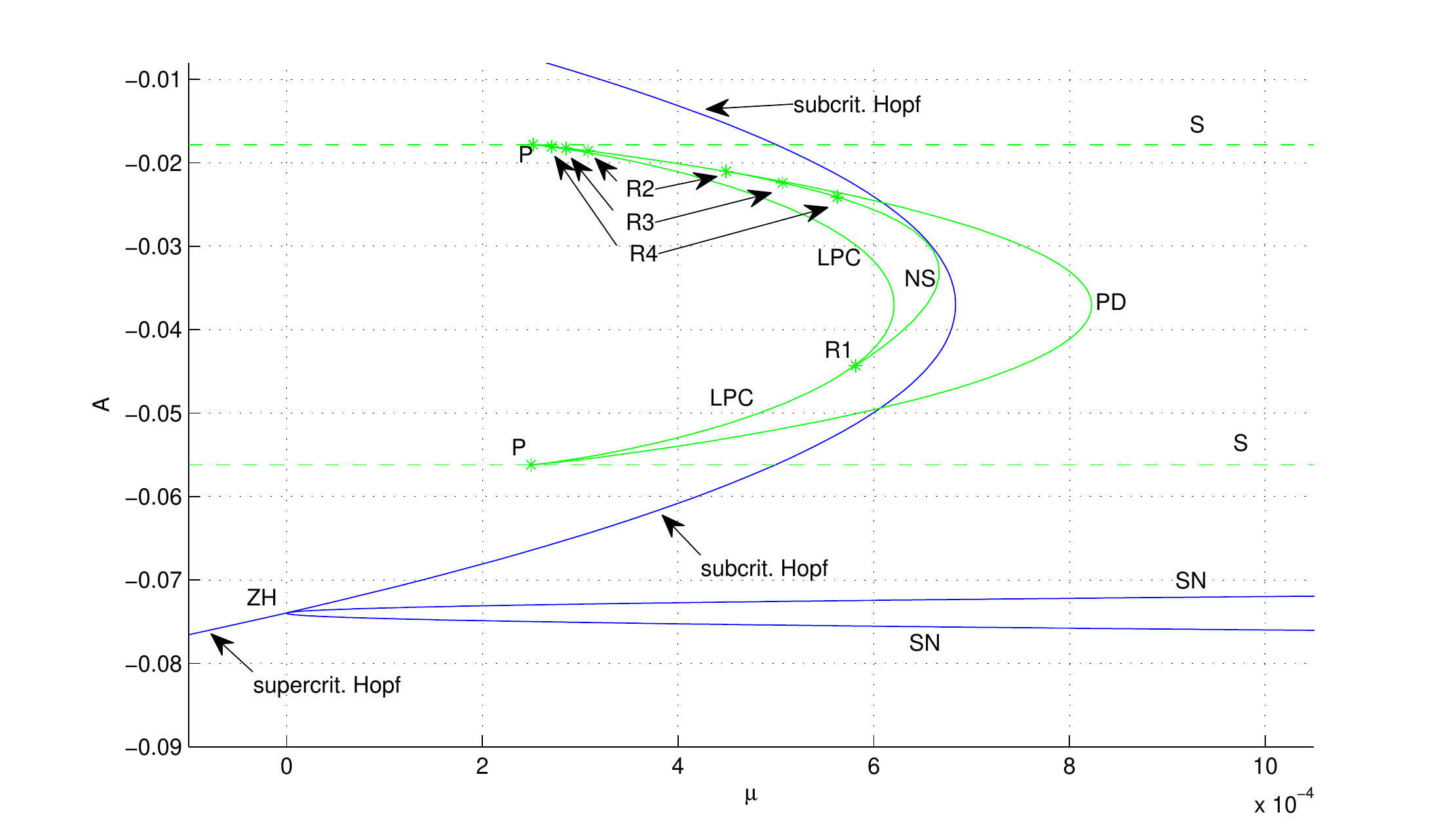} 
\caption{Region Va/VIa: $(\mu,A)$ bifurcation diagram for $(B,C)=(0.001,0.074)$. Relative to figure ~\ref{fig:bdiagB001C077_all_aexam}, the resonances have crossed the Hopf curve, and a second torus bifurcation curve with R2, R3 and R4 resonances appeared. }
\label{fig:bdiagB001C074_all_aexam}
\end{center}
\end{figure}

\begin{figure}[hptb]
\begin{center}
\includegraphics[width=1.0\textwidth]{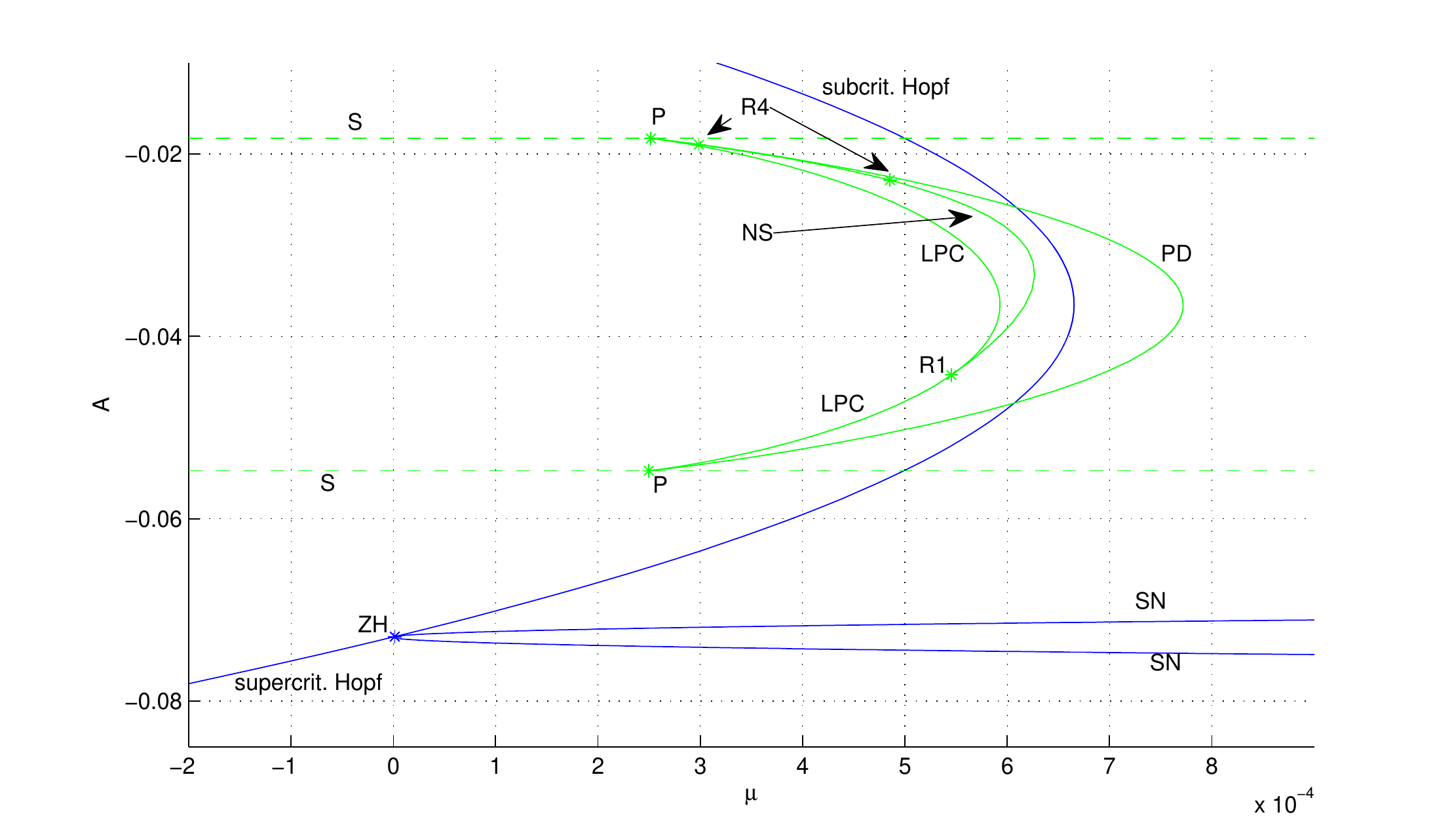} 
\caption{Region VIa: $(\mu,A)$ bifurcation diagram for $(B,C)=(0.001,0.073)$. Although the torus bifurcation curve and period-doubling curve are at times very close, we did not detect intersections of the two curves.  }
\label{fig:NEW_bdiag_B001C073_param27}
\end{center}
\end{figure}

\begin{figure}[hptb]
\begin{center}
\includegraphics[width=1.0\textwidth]{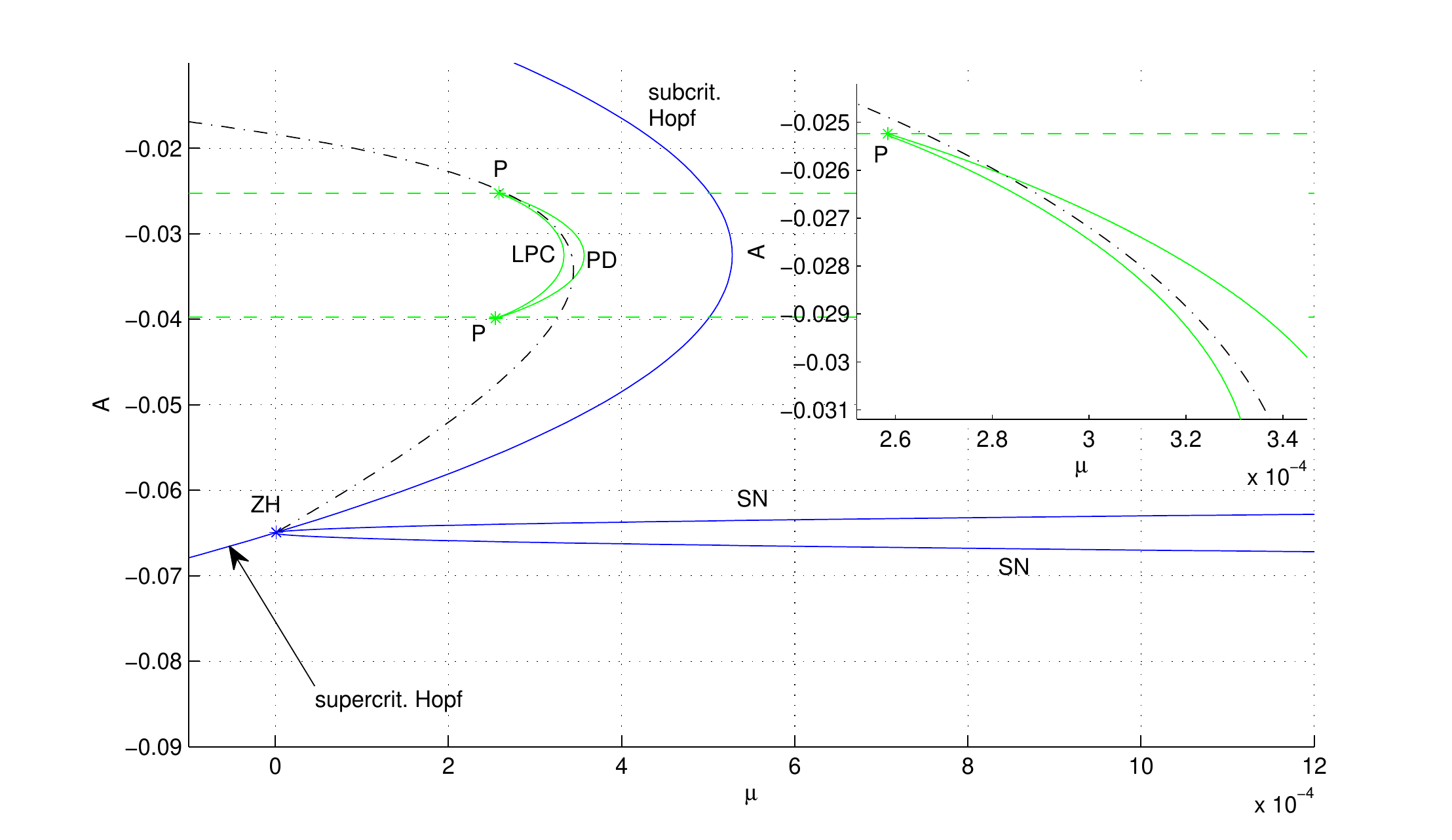} 
\caption{Region VIIa: $(\mu,A)$ bifurcation diagram for $(B,C)=(0.001,0.065)$. In this diagram, the curve $L_{(B,C)}$, explained in section 4 and not included in other $(\mu,A)$ bifurcation diagrams, is plotted with a dashed black line. The inset shows that the curve lies on the branch of the periodic orbit that has real multipliers. Note that while the singular Hopf periodic orbit does not undergo any torus bifurcations for these values of $B$ and $C$, there may be other $(\mu,A)$ bifurcation diagrams in region IVa with torus bifurcations and the same kinds of codimension 2 bifurcations of periodic orbits that occur in region IIIa. }
\label{fig:bdiagB001C065_all_aexam}
\end{center}
\end{figure}

\begin{figure}[hptb]
\begin{center}
\includegraphics[width=1.0\textwidth]{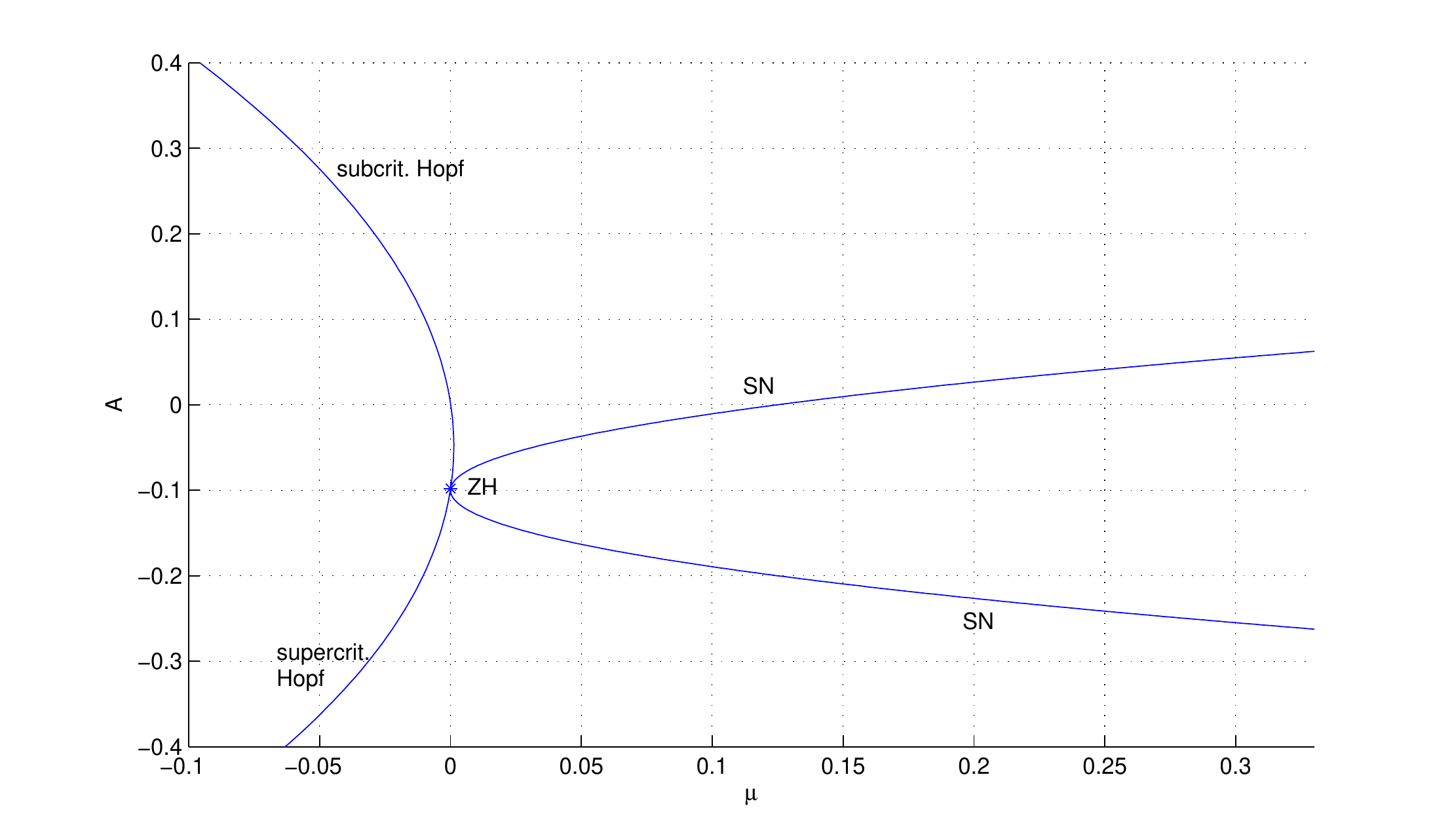} 
\caption{Region VIIIa: $(\mu,A)$ bifurcation diagram for $(B,C)=(0.02,0.1)$.}
\label{fig:bdiagB02C1_all_aexam}
\end{center}
\end{figure}

\clearpage

\section{Locations of non-transverse two parameter dimensional slices}

\begin{figure}[hptb]
\begin{center}
\includegraphics[width=1.0\textwidth]{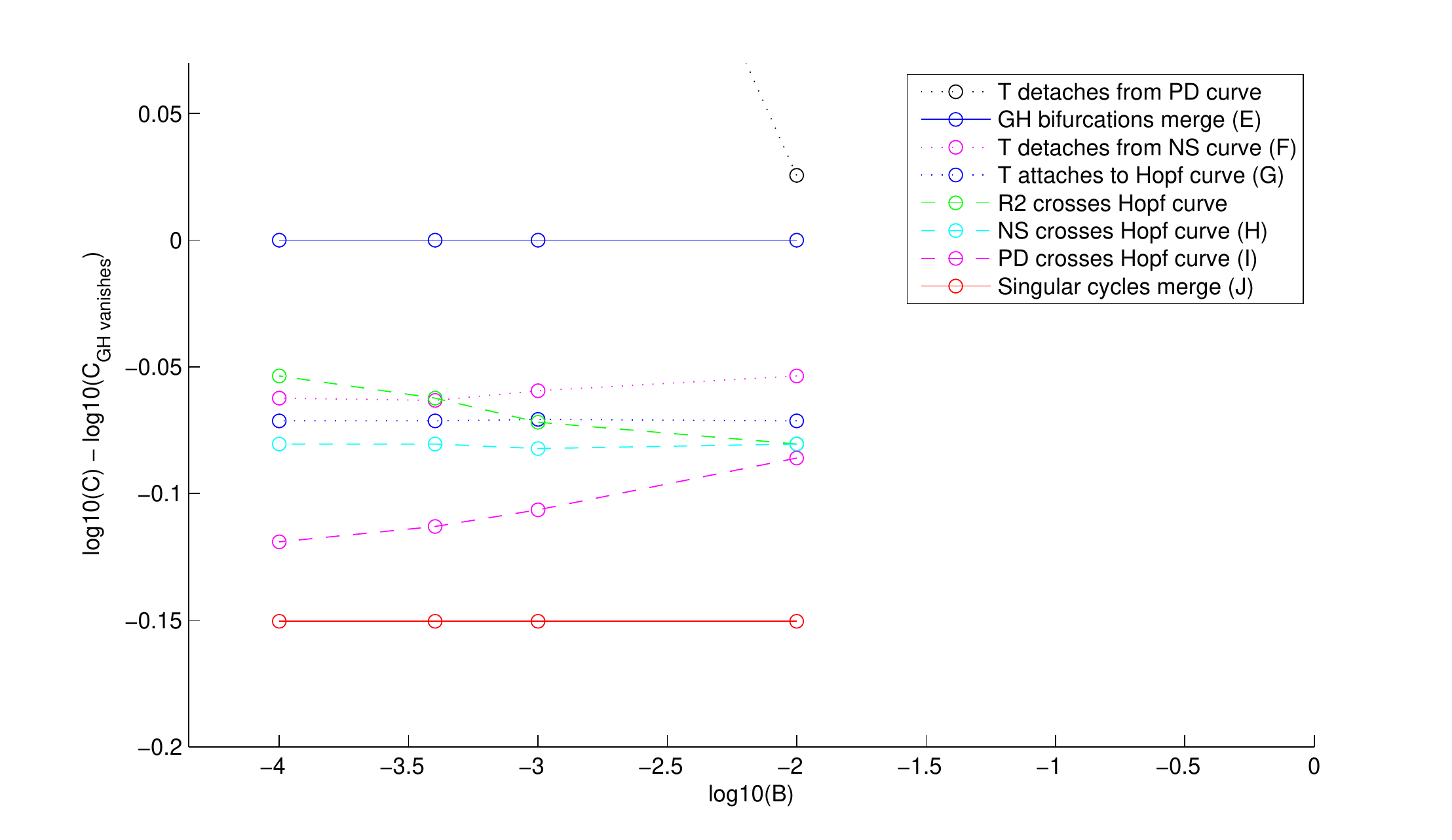} 

 \caption{Approximate positions of degeneracies of $(\mu,A)$ bifurcation diagrams in $(B,C)$ parameter space for $B,C>0$. \eJG{The legend lists eight types of degeneracies where $(\mu,A)$ slices of the parameter space fail to be transverse to codimension two bifurcation surfaces.
}
}
\label{fig:codim3_events}
\end{center}
\end{figure}

Inspecting $(\mu,A)$ bifurcation diagrams for a large number of parameters $(B,C)$, it is easy to see that the catalog of diagrams \eJGold{completely describes the bifurcations of the types considered in this paper} for $B<0$. The completeness of the diagrams for $B>0$ was verified by computing curves \eJGold{that separate the $(B,C)$ plane into regions with qualitatively different
$(\mu,A)$ bifurcation diagrams as shown in figure ~\ref{fig:codim3_events} .} \ePM{Where applicable, the curve labels introduced in section 4 are included in parentheses. }
For numerical reasons, $|B|$ was restricted to lie between $10^{-4}$ and $10^{-2}$. The positions of the folds of the generalized Hopf bifurcations and of the points where the slow flow has singular cycles are calculated using asymptotic relations. All other data points were obtained from sequences of $(\mu,A)$ \eJGold{bifurcation diagrams}. The error in the ordinates of the data points is significant, but the order of the data points on the ordinate is correct for each value of $B$. Note how the dashed green line intersects several other curves, corresponding to \eJGold{further degeneracy in the set of $(\mu,A)$ bifurcation diagrams.}\\

Figure ~\ref{fig:bdiagBm01C25_all_aexam} illustrates the degeneracies that can occur in $(\mu,A)$ bifurcation diagrams if $C$ is too large relative to $|B|^{1/2}$. The saddle-node bifurcation is close to intersecting the dashed curve where the slow-flow has singular cycles. 

\newpage 
\section{Numerical methods for the computation of the tangency curve}

\ePMnew{This appendix describes the two numerical methods used for the computation of the tangency curve. The first method is a shooting method, used by Desroches et al  (\cite{mmo_paper}) to compute the tangency curve for the Koper model discussed in section 5. The second method uses collocation and the software package AUTO \cite{AUTO}. }

The tangency bifurcation marks the onset of a bistability of trajectories in the unstable manifold $W^u(E_f)$ of the saddle-focus $E_f$ near the origin: the repelling slow manifold $S_r$ separates trajectories that flow to $X=-\infty$ from those that remain in the fold region.
This observation motivates the use of a shooting algorithm to compute the position of the tangency curve in parameter space: given a set of parameters, a grid of initial conditions in a linear approximation of $W^u(E_f)$ is integrated numerically for a long time interval. 
The tangency curve separates the parameters where at least one trajectory reaches $X=-\infty$ from the parameters where all trajectories approach a bounded attractor. This curve can be computed by a 2-parameter predictor-corrector continuation method starting from an initial parameter on the tangency curve, where the correction may be implemented using interval bisection. 

The main disadvantage of this method is that it requires choices regarding the number of grid points and, more importantly, the length of the integration-time. Unless convergence to an attractor such as a limit cycle can be determined easily, this method cannot determine in an automated fashion whether a trajectory that has stayed in the fold region for a long period of time will remain in the fold region forever or leave it after a phase of transient behavior. For example, if $E_f$ is close to a Hopf bifurcation curve, and the real part of its eigenvalues are thus close to zero, convergence to the periodic orbit is slow and takes a very long time. In this case, the algorithm may give incorrect results unless the integration time is chosen adaptively. 

The run-time of the algorithm described above is  proportional to the number of grid points and the integration-time, as the bulk of the computation time of the algorithm is spent on numerical time-integration of trajectories that remain in the fold region for the entire integration time. Slight speed-ups can be achieved by choosing the order of time integrations for grids of initial conditions to start with points that are far apart in the fundamental domain. Sometimes, this locates a trajectory that escapes the fold region more quickly than choosing initial points consecutively along a segment of a fundamental domain of the linearized unstable manifold. We found that using more than 10 grid points did not improve the accuracy of the method significantly. Larger speed-ups are possible by parallelizing the algorithm. Since the tangency curve is smooth, parallelism can be implemented using a two-pass approach: in the first pass, a coarse approximation of the position of the curve is computed on a single processor with the algorithm described above, using large step-sizes, few grid points, and thus little computation time. The second pass interpolates between data points of the first pass, using larger numbers of interpolated points, more grid points and long integration times on multiple processors. This second pass is embarrassingly parallel. 
\\

The tangency curve can also be computed via continuation of a boundary value problem in AUTO. The boundary conditions require that trajectories begin on a fixed ray in the linear approximation of the equilibrium's unstable manifold and end on the parabola $Y=X^2+5$. An initial trajectory segment satisfying these requirements, obtained for example using the first method, can then be continued in one parameter. AUTO's fold detection determines when the continued solution folds, i.e.  $W^u(E_f)$ and $S_r$ intersect tangentially. Switching to a continuation of the fold, the tangency can now be continued in two parameters. Computing the tangency curve in AUTO is in general much faster than the shooting method, but the shooting method provides additional insight into the fate of the trajectories in \ePMnew{$W^u(E_f)$} before and after the tangency.  As with the shooting method, the boundary value algorithm breaks down close to the Hopf curve, unless the number of mesh points is increased appropriately. \\

\end{document}